\documentclass[reqno,a4paper,12pt]{amsart}

\numberwithin{equation}{section}
\usepackage{amsmath, amsthm, amssymb, amscd, accents, bm, amsfonts}
\usepackage{mathtools}
\usepackage{url}
\usepackage{mathrsfs,dsfont}
\usepackage{mathtools}
\usepackage{mdwlist}
\usepackage{paralist}
\usepackage{datetime}
\usepackage{hyperref}
\usepackage{colonequals}
\usepackage{soul}
\usepackage[dvipsnames]{xcolor}
\usepackage{enumitem}

\usepackage[sort,nocompress]{cite}
\newtheorem{thm}{Theorem}
\newtheorem{theorem}[thm]{Theorem}
\newtheorem*{theorem*}{Theorem}

\newtheorem{lemma}[thm]{Lemma}
\newtheorem{cor}[thm]{Corollary}
\newtheorem{prop}[thm]{Proposition}

\newtheorem{question}[thm]{Question}

\theoremstyle{definition}
\newtheorem{definition}[thm]{Definition}
\newtheorem{remark}[thm]{Remark}

\newtheorem{example}[thm]{Example}

\numberwithin{equation}{section}
\numberwithin{thm}{section}

\def\F{{\mathbb F}}
\def\N{{\mathbb N}}
\def\R{{\mathbb R}}

\def\C{{\mathbb C}}

\newcommand{\sgn}{{\rm sgn}\,}
\newcommand{\lspan}{{\mathrm{span}}}

\newcommand{\re}{{\mathrm{Re} \,}}

\newcommand{\im}{{\mathrm{Im} \,}}

\makeatletter
\newcommand\thankssymb[1]{\textsuperscript{\@fnsymbol{#1}}}
\allowdisplaybreaks

\begin{document}

\title[Stable Phase Retrieval from Operator Measurements]{Stable Phase Retrieval from Operator Measurements}

\author[D. Needell]{Deanna Needell}
\address{Department of Mathematics, University of British Columbia, Vancouver, BC, Canada, V6T 1Z2 \& Department of Mathematics, University of California, Los Angeles, 90095 USA}
\email{needell@math.ubc.ca}
\author[Y. Plan]{Yaniv Plan}
\address{Department of Mathematics, University of British Columbia, Vancouver, BC, Canada, V6T 1Z2}
\email{yaniv@math.ubc.ca}
\author[T. Szczepanski]{Tomasz Szczepanski}
\address{Department of Mathematics, University of British Columbia, Vancouver, BC, Canada, V6T 1Z2}
\email{tomasz@math.ubc.ca}
\author[Ö. Y{\i}lmaz]{Özgür Y{\i}lmaz}
\address{Department of Mathematics, University of British Columbia, Vancouver, BC, Canada, V6T 1Z2 \& PIMS-CNRS International Research Laboratory}
\email{oyilmaz@math.ubc.ca}

\date{\today}

\begin{abstract}
A known result of Cahill, Casazza and Daubechies [Trans. Amer. Math. Soc. Ser. B\textbf{3}(2016), 63-76] states that no frame in an infinite-dimensional Hilbert space does stable phase retrieval. We show that stable recovery becomes possible for more general operator measurements, but only at a single exponent. Given a system of bounded operators $\{A_j\}_{j\in J}$ on a real or complex Hilbert space $H$, we study recovery of $f\in H$ up to a global phase from the quadratic measurements $\langle A_j f,f\rangle$ and from their phase preserving square roots.

We consider $p$-Bessel systems ($2\leq p<\infty$) for which we measure the quadratic data in $\ell_{p/2}$ and the square root data in $\ell_p$. In infinite dimensions, we show that whenever $p\neq 4$, no $p$-Bessel system of operators does stable phase retrieval. In the case of $p=4$, we show that stability still fails on any system of semidefinite operators. Nevertheless, we construct a system of rank-one operators, not all self-adjoint, for which stable recovery is possible. We also show that stability can be achieved on a system of self-adjoint operators of rank-two.

For both types of measurements, we show that the stability can be tested on orthogonal pairs. This reduction preserves the optimal constant for quadratic measurements and, when the operators are semidefinite, for their phase preserving square roots. We deduce that in finite dimensions, uniqueness of recovery implies its stability. 

Finally, we show that whenever only the moduli of the measurements are retained, both the orthogonal reduction and the implication from uniqueness to stability can fail, already in dimension two. We also give second moment measurements in cryo-electron microscopy (cryo-EM) as a motivation for the quadratic model.
\end{abstract}

\maketitle

\section{Introduction}
Throughout, let $H$ be a separable Hilbert space with inner product $\langle \cdot, \cdot \rangle$ over a real or complex scalar field $\mathbb{F}$. Let $J$ be a countable index set (finite or infinite). Recall, we say $\{g_j\}_{j\in J}\subseteq H$ is a \textit{frame} if there exist \textit{frame bounds} $0<A\leq B<\infty$ such that
    $$A\|f\|^2\leq \sum_{j\in J}|\langle f,g_j\rangle|^2\leq B\|f\|^2, \qquad f\in H.$$
While frames allow for stable recovery of $f\in H$ from linear measurements $(\langle f,g_j\rangle)_{j\in J}$ (see, e.g., \cite[Section 1.8]{casazza2013} for an overview of applications of frames in quantization, compressed sensing, and other areas), we consider a model where only the information of the magnitudes $|\langle f,g_j\rangle|$ is available. The phase retrieval problem asks what additional assumptions on a frame are needed to recover $f$ upon such loss of phase data. Prominent examples where such situation occurs include X-ray crystallography \cite{elser2003solution}, diffraction imaging \cite{miao2015beyond,shechtman2015phase,thibault2008high}, and quantum mechanics \cite{orl1994phase,raymer1997measuring} or signal processing and data analysis \cite{mallat2016understanding}.

For a frame $\{g_j\}_{j\in J}\subseteq H$ define $\mathcal{A}:H\to \ell_2(J)$ by $\mathcal{A}(f)=(|\langle f,g_j\rangle|)_{j\in J}.$ Note that the existence of an upper frame bound of $\{g_j\}_{j\in J}$ guarantees $\mathcal{A}(f)\in \ell_2(J)$ for every $f\in H$. 

\begin{definition}
    We say a frame $\{g_j\}_{j\in J}\subseteq H$ does \textit{phase retrieval} if for every $f,h\in H$,
        $$\mathcal{A}(f)=\mathcal{A}(h) \implies f=\theta h, \text{ for some } |\theta|=1.$$
    Furthermore, we say it does \textit{stable phase retrieval} (SPR) with constant $C>0$ if
        $$\min_{|\theta|=1}\|f-\theta h\|\leq C\|\mathcal{A}(f)-\mathcal{A}(h)\|_2, \qquad f,h\in H.$$
    \end{definition}

Uniqueness and stability in phase retrieval have been studied using methods from a wide range of viewpoints, for instance, the perspective of functional analysis \cite{Freeman2024,christ2022examples,freeman2025cahill,aldroubi2020phaseless,alharbi2023stable}, convex geometry \cite{bianchi2011phase}, complex analysis \cite{grohs2021stable,jaming2014uniqueness,wellershoff2024phase}, sampling theory \cite{grohs2025phase,jaming2024gabor,ramos,alaifari2024unique} or abstract algebra \cite{bendory2022algebraic}. We also refer to the survey  \cite{grohskopp2020phase}.

In this paper, we study a generalized model of phase retrieval, when frame measurements are replaced with operator measurements. Since
$$|\langle f,g_j\rangle|^2 = \langle g_jg_j^*f,f\rangle =: \langle A_jf,f\rangle,$$
it is natural to replace positive rank-one $A_j=g_jg_j^*$ with a more general operator. This motivates recovery from quadratic measurements (or their signed square roots) given by $\langle A_jf ,f\rangle$, 
for a family $\{A_j\}_{j\in J}$ from the space $B(H)$ of bounded linear operators on $H$.

\subsection{Phase retrieval via operator measurements}

Let $\{A_j\}_{j\in J}\subseteq B(H)$. Our goal is to find suitable families of operators that allow for stable recovery of a signal $f\in H$ from quadratic measurements $\langle A_jf,f\rangle$, or their signed square roots. In the case of frame measurements, the well-definedness of the measurement map $\mathcal{A}$ comes from the existence of an upper frame bound, also known as a \textit{Bessel bound}. It is natural to impose an analogous condition on the system $\{A_j\}_{j\in J}$, that generalizes this notion. More precisely, suppose $\{g_j\}_{j\in J}\subseteq H$ has a Bessel bound $0<B<\infty$, that is, 
$$\sum_{j\in J}|\langle f,g_j\rangle|^2\leq B\|f\|^2, \qquad f\in H.$$
By setting $A_j=g_jg_j^*$, the above is equivalent to
$$\sum_{j\in J}|\langle A_jf,f\rangle|\leq B\|f\|^2, \qquad f\in H.$$
Therefore, given $p\in [2,\infty)$, we say $\{A_j\}_{j\in J}\subseteq B(H)$ is a \textit{$p$-Bessel system} (\textit{$p$-Bessel}, for short) if there exists $0<B<\infty$ such that 
    $$\sum_{j\in J}|\langle A_jf,f\rangle|^{p/2}\leq B\|f\|^p, \qquad f\in H.$$

We now introduce the two main maps investigated in this paper. 

\begin{definition}
    Let $\{A_j\}_{j\in J}\subseteq B(H)$ be $p$-Bessel. For every $j$ define
    \begin{enumerate}
        \item $\beta_j:H\to \F$ by
        $$\beta_j(f):=\langle A_jf,f\rangle,$$
        \item $\alpha_j:H\to \F$ by 
    $$
        \alpha_j(f):=
        \begin{cases}
        \dfrac{\langle A_jf,f\rangle}
        {\sqrt{|\langle A_jf,f\rangle|}},
        & \text{if } \langle A_jf,f\rangle\neq 0,\\[2mm]
        0,
        & \text{if } \langle A_jf,f\rangle=0.
        \end{cases}
    $$
    \end{enumerate}
    Moreover, define the measurement maps
    \begin{enumerate}
    \item $\beta:H\to \ell_{p/2}(J)$ by
        $$\beta(f)=(\beta_j(f))_{j\in J}, \qquad f\in H,$$
    \item $\alpha:H\to \ell_p(J)$ by 
    $$\alpha(f) = (\alpha_j(f))_{j\in J}, \qquad f\in H.$$
    \end{enumerate}
\end{definition}

Note that whenever $H$ is a real Hilbert space or $A_j$ is self-adjoint, $\alpha_j$ can be equivalently defined via
    $$\alpha_j(f)=\operatorname{sgn}\bigl(\langle A_jf,f\rangle\bigr)\sqrt{|\langle A_jf,f\rangle|},$$
where we use the convention that $\sgn(0)=0$.

The map $\alpha$ directly generalizes the classical phase retrieval measurement map $\mathcal{A}$ in the sense that whenever $\{A_j\}_{j\in J}\subseteq B(H)$ is given by $A_j=g_jg_j^*$ for some frame $\{g_j\}_{j\in J}\subseteq H$, we have $\alpha(f)=\mathcal{A}(f)$ for all $f\in H$. While our primary interest is in recovery through $\alpha$, we also introduce $\beta$ because of its linearity in $ff^*$, which simplifies the analysis, as well as its connection with the cryo-electron microscopy model described in the next section. As we will show, stability estimates for recovery through $\beta$ imply corresponding stability estimates for recovery through $\alpha$.

Using the maps $\alpha$ and $\beta$, we generalize the notion of phase retrieval.

\begin{definition} 
\label{def:spr}
    Let $\{A_j\}_{j\in J}\subseteq B(H)$ be a $p$-Bessel system.
    \begin{enumerate}
        \item We say that $\{A_j\}_{j\in J}$ does \textit{$\beta$-phase retrieval} ($\beta$-PR) if for every $f,h\in H$,
        $$\beta(f)=\beta(h) \implies f=\theta h, \text{ for some } |\theta|=1.$$
    Moreover, we say it does $(p,\beta)$-stable phase retrieval ($(p,\beta)$-SPR) with constant $C>0$ if 
        $$\min_{|\theta|=1}\|f-\theta h\|\max_{|\theta|=1}\|f-\theta h\|\leq C\|\beta(f)-\beta(h)\|_{p/2}, \qquad f,h\in H.$$
        \item We say $\{A_j\}_{j\in J}$ does \textit{$\alpha$-phase retrieval} ($\alpha$-PR) if for every $f,h\in H$,
        $$\alpha(f)=\alpha(h) \implies f=\theta h, \text{ for some } |\theta|=1.$$
    Moreover, we say it does $(p,\alpha)$-stable phase retrieval ($(p,\alpha)$-SPR) with constant $C>0$ if 
        $$\min_{|\theta|=1}\|f-\theta h\|\leq C\|\alpha(f)-\alpha(h)\|_p,\qquad f,h\in H.$$
    \end{enumerate}
    We also simply say that $\{A_j\}_{j\in J}$ does $(p,\beta)$-SPR, if the $(p,\beta)$-SPR inequality holds for some $C>0$, and similarly for $(p,\alpha)$-SPR.        
\end{definition}

While it is easy to deduce that the notions of $\alpha$-PR and $\beta$-PR are equivalent (see \eqref{eq:alpha_PR_iff_beta_PR}), the relation between $(p,\alpha)$-SPR and $(p,\beta)$-SPR will be an important part of our investigation in Sections \ref{sec:inf_dim_stability} and \ref{sec:orthogonal_reductions_and_f_dim_stability}.

We also note that the choice of the min-max product corresponding to the map $\beta$ is motivated by the fact, that for every $z,t\in \F$ we have 
$$\min_{|\theta|=1}|z-\theta t|\max_{|\theta|=1}|z-\theta t| = \big||z|-|t|\big|\cdot \big(|z|+|t|\big) = \big||z|^2 - |t|^2\big|.$$
More generally for any $f,h\in H$, the min-max expression equals the trace norm distance
    $$\min_{|\theta|=1}\|f-\theta h\|\max_{|\theta|=1}\|f-\theta h\|=\|ff^*-hh^*\|_{\mathcal S_1},$$
which is the metric $d_1$ considered by Balan and Zou \cite[Section 2.1]{balan2016lipschitz}. In the real case, it reduces to $\|f-h\|\|f+h\|$, as used in the stability framework of Eldar and Mendelson \cite{eldar2014phase}.

Since for systems $A_j=g_jg_j^*$ we have $\alpha=\mathcal{A}$, the study of the classical phase retrieval problem corresponds to the study of the map $\alpha$ for positive rank-one operators associated with a frame.
The uniqueness problem for the map $\alpha$ has also been investigated in the context of orthogonal projections on $\F^d$ \cite{Cahill2016}, positive operators on $\R^d$ \cite{Zhuang2019}, positive operators on $\C^d$ \cite{Balan2021, zhuang2019complex} and positive operators on a general Hilbert space \cite{Han2019}. On the other hand, the uniqueness problem for $\beta$ was studied in the setting $\R^d$ \cite{Zhuang2019}, $\F^d$ \cite{Wang2019}, self-adjoint operators on $\C^d$ \cite{Balan2021,zhuang2019complex} and operators on a general Hilbert space \cite{Han2019}. 

In the direction of stability, for finite families of real symmetric matrices, Zhuang \cite{Zhuang2019} established quantitative bi-Lipschitz estimates for $\beta$ and showed that any process that does $\beta$-PR automatically satisfies such estimate. Analogous result was established for positive semidefinite matrices for $\alpha$ in \cite{zhuang2019complex}. Balan and Dock \cite{Balan2021} studied local and global Lipschitz properties in the more general setting of matrix recovery via matrix measurements. They established global stability of $\beta$ under injectivity and showed that, for matrix signals with more than one column, analogous global lower Lipschitz bounds for $\alpha$ fail with respect to natural quotient metrics.

Our main goal is to investigate when operator measurements allow for stable recovery in infinite-dimensional Hilbert spaces. To motivate the infinite-dimensional setting, let $\{A_j\}_{j\in J}\subseteq B(H)$ be $p$-Bessel, and let $(H_n)_{n\in \N}$ be an increasing sequence of finite-dimensional subspaces with dense union in $H$. At resolution $n$, suppose an unknown function is represented by $f_n\in H_n$, and the available measurements satisfy
\begin{align}
\label{eq:noisy_measurements}
    \widehat b=\bigl(\langle A_jf_n,f_n\rangle\bigr)_{j\in J}+\eta,
\end{align}
where $\eta\in\ell_{p/2}(J)$ is the measurement error. Assume now that the measurement map restricted to $H_n$ satisfies $(p,\beta)$-SPR with constant $C_n>0$. Then for any reconstruction $h\in H_n$ we have
\begin{align*}
    \min_{|\theta|=1}\|f_n-\theta h\|\max_{|\theta|=1}\|f_n-\theta h\|&\leq C_n\|\widehat b-\eta-\beta(h)\|_{p/2}\\
    &\leq C_n\left(\|\eta\|_{p/2}+\|\widehat b-\beta(h)\|_{p/2}\right).
\end{align*}
Thus $C_n$ controls how measurement error and the data fitting residual translate into reconstruction error. In particular, if $h$ fits the observed data at least as well as $f_n$, then
    $$\min_{|\theta|=1}\|f_n-\theta h\|\max_{|\theta|=1}\|f_n-\theta h\|\leq 2C_n\|\eta\|_{p/2}.$$
As the resolution increases, the corresponding stability constants $C_n$ may grow, weakening the reconstruction error bound. We therefore ask whether a fixed measurement model can admit a stability constant independent of the discretization. Stability at each fixed resolution does not by itself ensure such a bound, whereas a stability inequality on $H$ holds with the same constant on every $H_n$. As an example, we next describe how moment based cryo-electron microscopy gives rise to noisy quadratic measurements of the form \eqref{eq:noisy_measurements}.

\subsection{Cryo-electron microscopy}

Phase retrieval through operator measurements arises naturally in moment based approaches to single particle cryo-electron microscopy (cryo-EM). The goal of cryo-EM is to reconstruct the 3D structure of a biological molecule from noisy 2D projection images $I_i$ of $N$ copies of the molecule observed at unknown orientations. The unknown structure is represented by its electrostatic potential, or by a finite-dimensional approximation $f\in\mathbb R^n$; see \cite{bendory2020single} for an overview.

We model each image as
\begin{equation*}
    I_i=B_{\xi_i}f+\varepsilon_i,\qquad i=1,\dots,N,
\end{equation*}
where $\varepsilon_i$ represents noise and $\xi_i$ records the particle's orientation and any other imaging parameters included in the model. The linear operator $B_{\xi_i}:\mathbb R^n\to\mathbb R^m$ combines rotation and projection, and may also incorporate translations and microscope filtering. Although the family $\xi\mapsto B_\xi$ is specified by the imaging model, the parameters $\xi_i$ are not observed.

One approach to reconstruction uses empirical moments of the images rather than estimating each particle's orientation. This approach originates in Kam's autocorrelation analysis \cite{kam1980reconstruction} and has been developed extensively through the method of moments \cite{sharon2020method,bendory2023autocorrelation}. To explain its connection with our framework, suppose that the pairs $(\xi_i,\varepsilon_i)$ are independent and identically distributed, that $\xi_i$ has a known distribution with $\mathbb E\|B_\xi\|^2<\infty$, and that the noise is centered, independent of $\xi_i$, and has known covariance $\Sigma=\mathbb E(\varepsilon_i\varepsilon_i^*)$.

For a fixed image-space vector $u\in\mathbb R^m$, the scalar $\langle I_i,u\rangle$ is observable; choosing $u=e_k$, for example, extracts the $k$-th pixel. Averaging these scalars gives first moment information, which can lose angular structure; see \cite[Section 2.3.1]{sharon2020method}. Second moments retain additional information through correlations between image values.

Under our assumptions, a generic image $I=B_\xi f+\varepsilon$ satisfies
    $$\mathbb E|\langle I,u\rangle|^2=\mathbb E_\xi|\langle B_\xi f,u\rangle|^2+\langle\Sigma u,u\rangle.$$
The signal contribution is a quadratic measurement:
\begin{align}
\label{eq:cryo-EM_positive}
    \mathbb E_\xi|\langle B_\xi f,u\rangle|^2=\langle A_uf,f\rangle,\qquad A_u:=\mathbb E_\xi(B_\xi^*uu^*B_\xi).
\end{align}
Subtracting the known noise contribution therefore gives the empirical measurement
    $$\widehat b_{N,u}:=\frac1N\sum_{i=1}^N|\langle I_i,u\rangle|^2-\langle\Sigma u,u\rangle=\langle A_uf,f\rangle+\eta_{N,u},$$
where $\eta_{N,u}\to0$ almost surely as $N\to\infty$. This error reflects both the finite sampling of imaging parameters and fluctuations involving image noise. The operator $A_u$ is positive semidefinite and determined by the imaging model and parameter distribution. Although each operator inside its expectation has rank at most one, the average can have higher rank.

Mixed second moments similarly yield measurements associated with self-adjoint operators that may be indefinite. For $u,v\in\mathbb R^m$, define
    $$A_{u,v}:=\frac12\mathbb E_\xi\!\left[B_\xi^*(uv^*+vu^*)B_\xi\right].$$
Then
    $$\mathbb E[\langle I,u\rangle\langle I,v\rangle]=\langle A_{u,v}f,f\rangle+\langle\Sigma u,v\rangle,$$
so the corrected empirical measurement is
    $$\widehat b_{N,u,v}:=\frac1N\sum_{i=1}^N\langle I_i,u\rangle\langle I_i,v\rangle-\langle\Sigma u,v\rangle=\langle A_{u,v}f,f\rangle+\eta_{N,u,v},$$
with $\eta_{N,u,v}\to0$ almost surely. Taking $u$ and $v$ to be coordinate vectors gives the entries of the empirical image second moment matrix after subtracting the noise covariance. Thus both positive semidefinite and potentially indefinite self-adjoint operators occur naturally in this setting.

Choosing finitely many ordinary or mixed moments leads to the noisy quadratic measurement model
\begin{align*}
    \widehat b_N=\bigl(\langle A_jf,f\rangle\bigr)_{j\in J}+\eta_N.
\end{align*}
as in \eqref{eq:noisy_measurements}. Therefore, whenever the operator family satisfies the corresponding $(p,\beta)$-SPR inequality, its stability constant controls how errors in the above noisy measurements translate into error in the recovered signal up to global phase.

Recovery from these quadratic measurements is not automatically unique. In Kam's formulation with uniformly distributed orientations, second moments determine certain coefficient matrices only up to orthogonal factors. Resolving these factors is the \emph{orthogonal matrix retrieval} problem studied in \cite{bhamre2015orthogonal}. These ambiguities generally exceed the unavoidable global orientation ambiguity of the molecular structure. Additional information, structural assumptions, or higher-order moments are therefore needed in this setting. Nonuniform orientation distributions and sparsity assumptions lead to further recovery results under suitable hypotheses \cite{sharon2020method,bendory2023autocorrelation}. Related stability results for recovery from Gram matrices under structural assumptions have also been established in the framework of generalized phase retrieval over compact groups \cite{amir2026}.

\subsection{Main results}

A known theorem of J. Cahill, P. Casazza and I. Daubechies shows that no frame in an infinite-dimensional Hilbert space can do stable phase retrieval \cite[Theorem 2.2]{cahill2016phase}. By contrast, we show that stable recovery in infinite-dimensional Hilbert spaces is possible for the maps $\alpha$ and $\beta$ and is attained on a system of indefinite $4$-Bessel operators. We also show that the existence of such a system is extremely restrictive on the exponent $p\in [2,\infty)$: whenever $p\neq 4$, we show that no $p$-Bessel system achieves stability. We also show that stability cannot be achieved on any $4$-Bessel system consisting exclusively of semidefinite operators. In the setting of finite and infinite-dimensional Hilbert spaces, we then investigate necessary and sufficient conditions under which the stability constant can be tested on the set of orthogonal vectors. As an application of this method, we show that in finite-dimensional Hilbert spaces, any $p$-Bessel system of bounded operators that does $\beta$-PR and $\alpha$-PR automatically does $(p,\beta)$-SPR and $(p,\alpha)$-SPR.

We now list the main statements proved in this work. As we will observe, the existence of a $p$-Bessel system that does $(p,\alpha)$-SPR, implies the existence of a real-linear embedding of an infinite-dimensional Hilbert space into $\ell_{p/2}$. Since such an embedding is impossible whenever $p\neq 4$, we get the following.
\begin{theorem}
\label{thm:alpha_beta_failure_on_strong_Bessel}
    Let $H$ be an infinite-dimensional Hilbert space and let $p\in [2,\infty)$ satisfy $p\neq 4$. Then no $p$-Bessel system does $(p,\alpha)$-SPR or $(p,\beta)$-SPR.
\end{theorem}

In particular, Theorem \ref{thm:alpha_beta_failure_on_strong_Bessel} gives the impossibility of retrieval via frames of Cahill, Casazza and Daubechies. Indeed, since every frame is associated with a system of $2$-Bessel operators, their result follows immediately.

The next three statements are devoted to the study of the case when $p=4$. When $p=4$, the existence of an embedding of an infinite-dimensional Hilbert space into $\ell_{p/2}$ is no longer an obstruction. The first question becomes if systems of semidefinite operators, such as those arising in \eqref{eq:cryo-EM_positive}, can perform stable recovery in infinite-dimensional Hilbert spaces. The following shows that the answer is negative.

\begin{theorem}
\label{thm:alpha_beta_semidefinite_failure}
    Let $H$ be an infinite-dimensional Hilbert space. Then no $4$-Bessel system consisting only of semidefinite operators does $(4,\alpha)$-SPR or $(4,\beta)$-SPR.
\end{theorem}

It becomes a natural question if stability can be achieved on a $4$-Bessel system, if we drop the assumptions that every operator is semidefinite. We show that it is possible for rank-one operators, assuming not all of them are self-adjoint. 

\begin{prop}
\label{prop:beta_alpha_stability_on_rank_1}
    Let $H$ be an infinite-dimensional Hilbert space. There exists a $4$-Bessel system of rank-one operators, not all self-adjoint, that does $(4,\beta)$-SPR with constant $\sqrt{2}$. Moreover, it does $(4,\alpha)$-SPR with constant $2$.
\end{prop}

To achieve stability on an indefinite system of self-adjoint operators, such system must also contain operators of rank-two or higher. By symmetrizing and perturbing operators from Proposition \ref{prop:beta_alpha_stability_on_rank_1}, we obtain the following.

\begin{theorem}
\label{thm:4_Bessel_stability}
    Let $H$ be an infinite-dimensional Hilbert space.
    \begin{enumerate}
        \item There exists a $4$-Bessel system of self-adjoint operators of rank at most two doing $(4,\beta)$-SPR with constant $\sqrt{2}$. Moreover, it also does $(4,\alpha)$-SPR with constant $2$.
        \item Let $r\geq 2$ be an integer. For every $\varepsilon\in (0,1)$ there exists a $4$-Bessel system of self-adjoint operators of rank $r$ doing $(4,\beta)$-SPR with constant $\frac{\sqrt{2}}{1-\varepsilon}$. Moreover, it also does $(4,\alpha)$-SPR with constant $2\frac{\sqrt{1+\varepsilon}}{1-\varepsilon}$.
    \end{enumerate}
\end{theorem}

Next, a known result in classical phase retrieval states that it suffices to verify the stability inequality on orthogonal pairs, with no change in the constant \cite[Theorem~2.1]{Alharbi2024}. We establish an analogous orthogonality reduction principle for the map $\beta$.

\begin{theorem}
\label{thm:beta_sufficient_on_ortho_pairs}
    Let $H$ be a Hilbert space and let $p\in [2,\infty)$. A $p$-Bessel system $\{A_j\}_{j\in J} \subseteq B(H)$ does $(p,\beta)$-SPR with constant $C>0$ if and only if the $\beta$-SPR inequality holds on orthogonal pairs with the same constant.
\end{theorem}

We also show that a similar reduction holds for the map $\alpha$, although in this case the resulting constant might be worse.

\begin{theorem}
\label{thm:alpha_sufficient_on_ortho_pairs}
    Let $H$ be a Hilbert space and let $p\in [2,\infty)$. Suppose $\{A_j\}_{j\in J} \subseteq B(H)$ is a $p$-Bessel system. 
    \begin{enumerate}
        \item $\{A_j\}_{j\in J}$ does $(p,\alpha)$-SPR with constant $C>0$ if and only if the $\alpha$-SPR inequality holds on orthogonal pairs for some constant.
        \item If, in addition, $\{A_j\}_{j\in J} \subseteq B(H)$ consists of semidefinite operators, $\{A_j\}_{j\in J}$ does $\alpha$-SPR with constant $C>0$ if and only if the $\alpha$-SPR inequality holds on orthogonal pairs with the same constant.
    \end{enumerate}
\end{theorem}

These reduction principles have two consequences. First, $\alpha$ and $\beta$ stability notions are equivalent, up to a constant, depending on the given system of operators.

\begin{cor}
\label{cor:beta_and_alpha_are_equivalent}
    Let $H$ be a Hilbert space and let $\{A_j\}_{j\in J}\subseteq B(H)$ be a $p$-Bessel system for $p\in [2,\infty)$. Then $\{A_j\}_{j\in J}$ does $(p,\alpha)$-SPR if and only if it does $(p,\beta)$-SPR.
\end{cor}

To motivate the second consequence of the orthogonality reduction principle, we note that in classical phase retrieval for frames, uniqueness of recovery implies stability in finite-dimensional Hilbert spaces \cite[Proposition~1.4]{cahill2016phase}. For finite systems of operators on finite-dimensional Hilbert spaces, the same implication was established for $\beta$ when the operators are self-adjoint, and for $\alpha$ when they are positive semidefinite \cite[Theorems~2.3-2.4]{Zhuang2019,zhuang2019complex}. Using the orthogonal reduction, we show that this is true for any $p$-Bessel system of operators. 

\begin{theorem}
\label{thm:beta_and_alpha_PR_in_finite_dim}
    Let $H$ be a finite-dimensional Hilbert space and let $\{A_j\}_{j\in J}\subseteq B(H)$ be a $p$-Bessel system for $p\in [2,\infty)$. Then $\{A_j\}_{j\in J}$ does $(p,\beta)$-SPR (and therefore $(p,\alpha)$-SPR) with some constant $C>0$ if and only if it does $\beta$-PR (and therefore $\alpha$-PR).
\end{theorem}

In the final section of this work we then investigate unsigned versions of the maps $\alpha$ and $\beta$, denoted by $|\alpha|$ and $|\beta|$, for which the information about the local phase is also lost. We show that the systems from Proposition \ref{prop:beta_alpha_stability_on_rank_1} and Theorem \ref{thm:4_Bessel_stability} do not achieve stable recovery through these maps. Furthermore, we show that even in finite-dimensional Hilbert spaces, an analogue of Theorem \ref{thm:beta_and_alpha_PR_in_finite_dim} does not hold in general, since uniqueness need no longer imply stability. Similarly, we show that the orthogonality reduction as described by Theorems \ref{thm:beta_sufficient_on_ortho_pairs} and \ref{thm:alpha_sufficient_on_ortho_pairs} can fail.

\subsection{Outline}
The organization of this paper is as follows. In Section \ref{sec:inf_dim_stability}, we investigate basic properties of the measurement maps and $p$-Bessel systems. This will allow us to prove, in the infinite-dimensional Hilbert setting, Theorems \ref{thm:alpha_beta_failure_on_strong_Bessel}, \ref{thm:alpha_beta_semidefinite_failure} together with Proposition \ref{prop:beta_alpha_stability_on_rank_1} and Theorem \ref{thm:4_Bessel_stability}.

In Section \ref{sec:orthogonal_reductions_and_f_dim_stability}, we prove the orthogonality reductions given by Theorems \ref{thm:beta_sufficient_on_ortho_pairs} and \ref{thm:alpha_sufficient_on_ortho_pairs}. We then use them to prove Corollary \ref{cor:beta_and_alpha_are_equivalent} and Theorem \ref{thm:beta_and_alpha_PR_in_finite_dim}.

Finally, in Section \ref{sec:unsigned_maps} we investigate whether stability can be achieved when considering the unsigned version of the measurement maps. We show that the systems of operators from Proposition \ref{prop:beta_alpha_stability_on_rank_1} and Theorem \ref{thm:4_Bessel_stability} that achieve signed stability fail the unsigned version of it. We also show how the orthogonality reduction principles can fail for these maps. Moreover, in finite-dimensional Hilbert spaces, stability need not be achieved for all systems doing phase retrieval.

\section{Infinite-dimensional stability}
\label{sec:inf_dim_stability}

Our main goal for this section is to provide necessary conditions for a system of operators to achieve stability in infinite-dimensional Hilbert spaces. We will then give examples of systems for which the stability inequality is satisfied. In this section $H$ stands for a Hilbert space over real or complex scalar field $\F$, finite or infinite-dimensional, unless specified. 

We start with a simple observation that finite measurement systems are always $p$-Bessel for every $p\in [2,\infty)$. Indeed, suppose $|J|<\infty$ and $\{A_j\}_{j\in J}\subseteq B(H)$. For every $f\in H$ we have
    $$\sum_{j\in J}|\langle A_jf,f\rangle|^{p/2} \leq \sum_{j\in J}\|A_j\|^{p/2}\|f\|^p \leq B\|f\|^p,$$
for any positive $B\geq \sum_{j\in J}\|A_j\|^{p/2}$.

Next, we show how the Bessel condition depends on $p$. 

\begin{lemma}
\label{lem:p_Bessel_relation}
    Let $2\leq p_1\leq p_2<\infty$. If $\{A_j\}_{j\in J}\subseteq B(H)$ is $p_1$-Bessel with constant $B$, then it is $p_2$-Bessel with constant $B^{p_2/p_1}$.
\end{lemma}
\begin{proof} 
    Suppose that $\{A_j\}_{j\in J}$ is $p_1$-Bessel with constant $B$. Fix $f\in H$. Since $\{A_j\}_{j\in J} $ is $p_1$-Bessel, we have
        $$\sum_{j\in J}|\langle A_jf,f\rangle|^{p_1/2}\leq B\|f\|^{p_1}.$$
    Since $\ell_{p_1/2}\subseteq \ell_{p_2/2} $, we have 
        $$\left(\sum_{j\in J} |\langle A_jf,f\rangle|^{p_2/2}\right)^{2/p_2} \leq \left(\sum_{j\in J} |\langle A_jf,f\rangle|^{p_1/2}\right)^{2/p_1}.$$ 
    Raising both sides to the power $p_2/2$ gives 
    \begin{align*}
        \sum_{j\in J} |\langle A_jf,f\rangle|^{p_2/2} &\leq \left( \sum_{j\in J} |\langle A_jf,f\rangle|^{p_1/2} \right)^{p_2/p_1} \leq (B\|f\|^{p_1})^{p_2/p_1} \\
        &= B^{p_2/p_1}\|f\|^{p_2}.
    \end{align*}
    Hence $\{A_j\}_{j\in J}$ is $p_2$-Bessel with constant $B^{p_2/p_1}$. 
\end{proof}

According to Lemma \ref{lem:p_Bessel_relation}, given $p_2\geq p_1$, every system that is $p_1$-Bessel is also $p_2$-Bessel. Nevertheless, since $\ell_{p_1}\subseteq \ell_{p_2}$, the stability relation is reversed: every such a system that does $(p_2,\alpha)$-SPR does $(p_1,\alpha)$-SPR (the same statement holds also for the map $\beta$). 

For a finite system $\{A_j\}_{j\in J}$ the situation is much more straightforward; every system is $p$-Bessel for every $p\in [2,\infty)$. Moreover, as the following lemma shows, if the stability holds for one exponent $p\in [2,\infty)$, it is satisfied for all of them.

\begin{lemma}
   Let $\{A_j\}_{j\in J}\subseteq B(H)$ be a system with $|J|<\infty$. Suppose there exists $p_0\in [2,\infty)$ such that $\{A_j\}_{j\in J}$ does $(p_0,\alpha)$-SPR with constant $C>0$. Then for every $p\in [2,\infty)$, $\{A_j\}_{j\in J}$ does $(p,\alpha)$-SPR with constant 
        $$C_p=C|J|^{\max\{0,\,1/p_0-1/p\}}.$$
\end{lemma}
\begin{proof}
    Fix $p\in [2,\infty)$ and $f\in H$. If $p\leq p_0$, then $\|f\|_{p_0}\leq \|f\|_{p}$. If $p> p_0$, applying Hölder’s inequality with conjugate exponents
        $$r=\frac p{p_0}, \qquad r'=\frac p{p-p_0}$$
    gives
    \begin{align*}
        \sum_{j\in J}|f_j|^{p_0}&\leq\left(\sum_{j\in J}|f_j|^{p_0r}\right)^{1/r}\left(\sum_{j\in J}1^{r'}\right)^{1/r'}\\
        &=\left(\sum_{j\in J}|f_j|^p\right)^{p_0/p}\cdot {|J|}^{1-p_0/p}.
    \end{align*}
    By taking the $p_0$-th root we obtain
        $$\|f\|_{p_0}\leq |J|^{1/p_0-1/p}\|f\|_p.$$
    Combining both cases we have that for every $p\in [2,\infty)$ and $f\in H$,
        $$\|f\|_{p_0}\leq |J|^{\max\{0,\,1/p_0-1/p\}}\|f\|_p.$$
    Now, suppose that $\{A_j\}_{j\in J}$ does $(p_0,\alpha)$-SPR with constant $C>0$. Fix $p\in [2,\infty)$. Then for every $f,h\in H$ we have
        $$\min_{|\theta|=1}\|f-\theta h\|\leq C\|\alpha(f)-\alpha(h)\|_{p_0}\leq C|J|^{\max\{0,\,1/p_0-1/p\}}\|\alpha(f)-\alpha(h)\|_p.$$
    Therefore, $\{A_j\}_{j\in J}$ does $(p,\alpha)$-SPR with constant $C_p=C|J|^{\max\{0,\,1/p_0-1/p\}}$.
\end{proof}

Suppose now that $\{A_j\}_{j\in J}\subseteq B(H)$ is $p$-Bessel for some $p\in [2,\infty)$. Observe that for every $j$ we have
$$|\alpha_j(f)|^p = |\beta_j(f)|^{p/2}.$$
This, together with 
$$\beta_j(f) = \alpha_j(f)|\alpha_j(f)|$$
allows us to deduce that 
\begin{align}
\label{eq:alpha_PR_iff_beta_PR}
    \alpha(f)=\alpha(h) \iff \beta(f)=\beta(h).
\end{align}

Therefore, $\{A_j\}_{j\in J}$ system does $\alpha$-PR if and only if it does $\beta$-PR. Our next goal is to show that stability for the map $\beta$ implies stability on $\alpha$ (in fact, we will see in Section \ref{sec:orthogonal_reductions_and_f_dim_stability} that these two notions are equivalent). For that purpose, we need the following.

\begin{lemma}
\label{lem:numbers_abs_value_inequality}
    Let $a,b\in \F$. Then 
    $$\Big| a|a|-b|b|\Big|\leq (|a|+|b|)|a-b|.$$
\end{lemma}
\begin{proof}
    Note that
    $$a|a|-b|b|=|a|(a-b)+b(|a|-|b|).$$
    Thus
    \begin{align*}
        \Big| a|a|-b|b|\Big| &= \Big| |a|(a-b)+b(|a|-|b|)\Big| \leq |a||a-b|+|b|\big||a|-|b|\big|\\
        &\leq |a||a-b|+|b||a-b| = (|a|+|b|)|a-b|.
    \end{align*}
\end{proof}

The second ingredient is that for every $f,h\in H$ we have
\begin{align}
\label{eq:max_phase_inequality}
    \|f\|+\|h\|\leq \sqrt{2}\max_{|\theta|=1}\|f-\theta h\|.
\end{align}
Indeed, note that $(\|f\|+\|h\|)^2 \leq 2(\|f\|^2+\|h\|^2).$ Hence 
    \begin{align*}
    \|f\|+\|h\|\leq \sqrt{2}\sqrt{\|f\|^2+\|h\|^2} &\leq \sqrt{2}\sqrt{\|f\|^2+\|h\|^2 + 2|\langle f,h\rangle|} \\
    &= \sqrt{2}\max_{|\theta|=1}\|f-\theta h\|.
    \end{align*}

The following statement shows the connection between $(p,\alpha)$-SPR and $(p,\beta)$-SPR.

\begin{prop}
\label{prop:beta_implies_alpha}
    Let $\{A_j\}_{j\in J} \subseteq B(H)$, and let $p\in [2,\infty)$. Suppose $\{A_j\}_{j\in J}$ is $p$-Bessel with constant $B$. If it does $(p,\beta)$-SPR with constant $C$, then it does $(p,\alpha)$-SPR with constant $C\sqrt{2}B^{1/p}.$ 
\end{prop}
\begin{proof}
    First, since $\{A_j\}_{j\in J}$ is $p$-Bessel, for every $f\in H$ we have
    $$\|\alpha(f)\|_p^p = \sum_{j\in J} |\langle A_jf,f\rangle|^{p/2}\leq B\|f\|^p$$
    thus $\|\alpha(f)\|_p \leq B^{1/p}\|f\|.$

    Next, observe that for every $j$
    $$\alpha_j(f)|\alpha_j(f)| = \beta_j(f).$$
    Hence, applying Lemma \ref{lem:numbers_abs_value_inequality} with $a=\alpha_j(f)$ and $b=\alpha_j(h)$ gives
    \begin{align*}
        |\beta_j(f)-\beta_j(h)| &= \Big| \alpha_j(f)|\alpha_j(f)| - \alpha_j(h)|\alpha_j(h)|\Big|\\
        &\leq (|\alpha_j(f)|+|\alpha_j(h)|)|\alpha_j(f)-\alpha_j(h)|.
    \end{align*}
    Taking the $\ell_{p/2}$ norm, applying Hölder's inequality with $\frac{2}{p}=\frac{1}{p}+\frac{1}{p}$ and applying  \eqref{eq:max_phase_inequality} gives
    \begin{align*}
        \|\beta(f)-\beta(h)\|_{p/2} &\leq \| (|\alpha_j(f)|+|\alpha_j(h)|)_{j\in J}\|_p\cdot \|\alpha(f)-\alpha(h)\|_p \\
        &\leq (\|\alpha(f)\|_p+\|\alpha(h)\|_p)\|\alpha(f)-\alpha(h)\|_p\\
        &\leq B^{1/p}(\|f\|+\|h\|)\|\alpha(f)-\alpha(h)\|_p\\
        &\leq \sqrt{2}B^{1/p}\max_{|\theta|=1}\|f-\theta h\|\|\alpha(f)-\alpha(h)\|_p
    \end{align*}
    Thus, if $\{A_j\}_{j\in J}$ does $(p,\beta)$-SPR with $C$, then 
    \begin{align*}
        \min_{|\theta|=1}\|f-\theta h\|\max_{|\theta|=1}\|f-\theta h\|&\leq C\|\beta(f)-\beta(h)\|_{p/2} \\
        &\leq C\sqrt{2}B^{1/p}\max_{|\theta|=1}\|f-\theta h\|\|\alpha(f)-\alpha(h)\|_p.
    \end{align*}
    If $\max_{|\theta|=1}\|f-\theta h\|=0$, then from  \eqref{eq:max_phase_inequality} we have $f=0=h$, so in particular $(p,\alpha)$-SPR inequality is trivially satisfied for this pair. If $\max_{|\theta|=1}\|f-\theta h\|\neq 0$, canceling this term on both sides gives 
    $$\min_{|\theta|=1}\|f-\theta h\|\leq C\sqrt{2}B^{1/p}\|\alpha(f)-\alpha(h)\|_p.$$
\end{proof}

Our next goal is to show that for $p\neq 4$ there are no $p$-Bessel systems that perform $(p,\alpha)$-SPR, and therefore $(p,\beta)$-SPR. For this purpose we start with an observation, that in infinite-dimensional Hilbert spaces, only infinite systems can do $\alpha$-PR and $\beta$-PR, and thus achieve stability for these maps.

\begin{lemma}
\label{lem:finite_systems_in_infinite_Hilbert}
    Let $H$ be infinite-dimensional and let $\{A_j\}_{j\in J}\subseteq B(H)$ satisfy $|J|<\infty$. Then $\{A_j\}_{j\in J}$ fails $\alpha$-PR and $\beta$-PR.
\end{lemma}
\begin{proof}
    Assume that $|J|=m<\infty$. Fix a unit vector $u\in H$. The space
        $$M=\operatorname{span}\{e,A_1u,\dots,A_mu,A_1^*u,\dots,A_m^*u\}$$
    is finite-dimensional. Since $H$ is infinite-dimensional, we can find a unit vector $v\in M^\perp$. For every $j$,
        $$\langle A_ju,v\rangle=0,\qquad \langle A_jv,u\rangle=\langle v,A_j^*u\rangle=0.$$
    Consequently,
        $$\langle A_j(u+v),u+v\rangle=\langle A_ju,u\rangle+\langle A_jv,v\rangle=\langle A_j(u-v),u-v\rangle.$$
    Thus
        $$\beta(u+v)=\beta(u-v),\qquad\alpha(u+v)=\alpha(u-v).$$
    Note that $u+v$ and $u-v$ are nonzero and orthogonal, since $u\perp v$ and $\|u\|=\|v\|=1$. Therefore $u+v\neq \theta(u-v)$ for all $|\theta|=1$.
\end{proof}

Next, to estimate the measurements of the map $\alpha$ using $\beta$, we need the following.

\begin{lemma}
\label{lem:inequality_for_alpha}
    For every non-zero $a,b\in \F$ we have 
        $$\Big|\frac{a}{\sqrt{|a|}}-\frac{b}{\sqrt{|b|}}\Big|^2\leq 2|a-b|.$$
    The statement also holds for $a=0$ or $b=0$, if we replace $\frac{a}{\sqrt{|a|}}$ or $\frac{b}{\sqrt{|b|}}$ by $0$, respectively.
\end{lemma}
\begin{proof}
    The cases $a=0$ or $b=0$ are immediate. Fix non-zero $a,b\in \F$. Let $x=\frac{a}{\sqrt{|a|}}$ and $y=\frac{b}{\sqrt{|b|}}$. Then $a=x|x|$ and $b=y|y|$. 
    Note that
        $$2\bigl(x|x|-y|y|\bigr)=(|x|+|y|)(x-y)+(|x|-|y|)(x+y)=:s+t.$$
    Since $(x-y)\overline{(x+y)}=|x|^2-|y|^2+2i\im(x\overline{y})$, we get $\re(s\overline{t})=(|x|^2-|y|^2)^2\geq 0$. Hence 
    \begin{align*}
      4\big|x|x|-y|y|\big|^2&=|s+t|^2=|s|^2+|t|^2+2\re (s\overline{t})\\
      &\geq |s|^2=(|x|+|y|)^2|x-y|^2,  
    \end{align*}
    or equivalently
    $$(|x|+|y|)|x-y|\leq 2\big|x|x|-y|y|\big|.$$
    Substituting back $a$ and $b$ gives
    $$(\sqrt{|a|}+\sqrt{|b|})\Big|\frac{a}{\sqrt{|a|}}-\frac{b}{\sqrt{|b|}}\Big|\leq 2|a-b|.$$
    We get
    \begin{align*}
    \Big|\frac{a}{\sqrt{|a|}}-\frac{b}{\sqrt{|b|}}\Big|^2 &= \Big|\frac{a}{\sqrt{|a|}}-\frac{b}{\sqrt{|b|}}\Big|\cdot \Big|\frac{a}{\sqrt{|a|}}-\frac{b}{\sqrt{|b|}}\Big| \\
    &\leq(\sqrt{|a|}+\sqrt{|b|})\Big|\frac{a}{\sqrt{|a|}}-\frac{b}{\sqrt{|b|}}\Big|\leq 2|a-b|.
    \end{align*}
\end{proof}

We are ready to prove Theorem \ref{thm:alpha_beta_failure_on_strong_Bessel}.
\begin{proof}[Proof of Theorem \ref{thm:alpha_beta_failure_on_strong_Bessel}]
    Let $H$ be infinite-dimensional, let $p\neq 4$ and let $\{A_j\}_{j\in J}\subseteq B(H)$ be $p$-Bessel with constant $B$. If $|J|<\infty$, $\{A_j\}_{j\in J}$ fails $\alpha$-PR and $\beta$-PR from Lemma \ref{lem:finite_systems_in_infinite_Hilbert}, therefore also $(p,\alpha)$-SPR and $(p,\beta)$-SPR for every $p$. Hence, from now on we assume that $|J|=\infty$
    
    It suffices to show that $\{A_j\}_{j\in J}$ fails $(p,\alpha)$-SPR, by Proposition \ref{prop:beta_implies_alpha}. Suppose, by contradiction, that $\{A_j\}_{j\in J}$ does $(p,\alpha)$-SPR with constant $C>0$. We will show that this implies existence of an infinite-dimensional Hilbert space that embeds into $\ell_{p/2}(J)$, which is impossible for $p\neq 4$.

    As a preliminary calculation, since for every $f\in H$ we have 
        $$\sum_{j\in J}|\langle A_jf,f\rangle|^{p/2}\leq B\|f\|^p,$$
    taking the $p/2$-th root gives 
    \begin{align}
    \label{eq:pth_root_for_Bessel}
        \|\beta(f)\|_{p/2}\leq B^{2/p}\|f\|^2.
    \end{align}
    Let $u\in H$ be a unit vector. We note that $u^\perp$ is an infinite-dimensional Hilbert space. Our goal is to show that $u^\perp$ embeds into $\ell_{p/2}(J)$ via a suitable operator $T:u^\perp \to \ell_{p/2}(J)$. To define this map, first observe that for every $v\in u^\perp$ we have
        \begin{align*}
            \langle A_j(u+v),u+v\rangle - \langle A_j(u-v),u-v\rangle &= 2(\langle A_ju,v\rangle+\langle A_jv,u\rangle).
        \end{align*}
    Therefore, for every $v\in u^\perp$,
        $$\frac{1}{2}\big(\beta(u+v)-\beta(u-v)\big)=(\langle A_ju,v\rangle+\langle A_jv,u\rangle)_{j\in J}.$$
    In particular, the map 
        $$T:u^\perp\to \ell_{p/2}(J), \qquad Tv = (\langle A_ju,v\rangle + \langle A_jv,u\rangle)_{j\in J}$$
    is well-defined. We note that $T$ is real-linear.
    
    We claim that for every unit $v\in u^\perp$, we have
    \begin{align}
    \label{eq:embedding_bounds}
        \frac{1}{C^2}\|v\|\leq \|Tv\|_{p/2}\leq 2B^{2/p}\|v\|.
    \end{align}
    Fix unit $v\in u^\perp$. Using \eqref{eq:pth_root_for_Bessel} we get 
    \begin{align*}
        2\|Tv\|_{p/2} &\leq \|\beta(u+v)\|_{p/2}+\|\beta(u-v)\|_{p/2}\leq B^{2/p}\|u+v\|^2 + B^{2/p}\|u-v\|^2\\
        &=2B^{2/p}(2\|u\|^2+2\|v\|^2)=4B^{2/p},
    \end{align*}
    which gives the upper bound in \eqref{eq:embedding_bounds}. Next, using the fact that $v\perp u$ and $\|v\|=1=\|u\|$, we have $u-v \perp u+v$, which implies 
        $$\min_{|\theta|=1}\|(u-v)-\theta(u+v)\|= 2.$$
    Since $\{A_j\}_{j\in J}$ does $(p,\alpha)$-SPR with constant $C>0$, we get 
        $$2\leq C\|\alpha(u+v)-\alpha(u-v)\|_p.$$
    Using Lemma \ref{lem:inequality_for_alpha}, for every $j$ we have
    \begin{align*}
        |\alpha_j(u+v)-\alpha_j(u-v)|^2 \leq 2|\beta_j(u+v)-\beta_j(u-v)|.
    \end{align*}
    Taking the $\ell_{p/2}$ norm gives 
        $$\|\alpha(u+v) - \alpha(u-v)\|_p^2 \leq 2\|\beta(u+v)-\beta(u-v)\|_{p/2}.$$
    Thus
    \begin{align*}
        4&\leq C^2\|\alpha(u+v)-\alpha(u-v)\|_p^2 \leq 2C^2\|\beta(u+v)-\beta(u-v)\|_{p/2}\\
        &=4C^2\|Tv\|_{p/2}.
    \end{align*}
    This proves the lower bound of \eqref{eq:embedding_bounds}.
    By real homogeneity of $T$, \eqref{eq:embedding_bounds} holds for all $v\in u^\perp$. 
    When $H$ is a real Hilbert space, this gives an embedding of an infinite-dimensional Hilbert space $u^\perp$ into $\ell_{p/2}(J;\R)$, which is impossible when $p\neq 4$, see \cite[Corollary 2.1.6]{albiac2006topics}. To arrive at the same contradiction when $H$ is a complex Hilbert space, we view $H$ as a real Hilbert space with the inner product $\langle x,y\rangle_\R = \re\langle x,y\rangle$. Consider the real-linear map $S:\ell_{p/2}(J;\C) \to \ell_{p/2}(J,\R)$ given by $Sz=(\re z, \im z)$. Using the equivalence of norms on $\R^2$, we can find constants $0<a\leq b<\infty$, such that 
    $$a\|z\|\leq \|Sz\|\leq b\|z\|, \qquad z\in \ell_{p/2}(J;\C).$$
    Therefore, since 
        $$\frac{a}{C^2}\|v\|\leq \|STv\|\leq 2B^{2/p}b\|v\|, \qquad v\in u^\perp,$$
    we embed $u^\perp$ into real $\ell_{p/2}(J,\R)$, giving a contradiction.
\end{proof}

When $p\neq 4$, Theorem \ref{thm:alpha_beta_failure_on_strong_Bessel} shows failure of stability when measured through any $p$-Bessel system. Next, we show that failure for semidefinite operators occurs also in the case $p=4$. We now provide a short argument why it suffices to check the failure for positive semidefinite operators. Suppose $\{A_j\}_{j\in J}\subseteq B(H)$ is a system consisting of semidefinite operators. Note that for every $j$ and $f\in H$ we have
\begin{align*}
    |\langle |A_j|f,f\rangle| = |\langle A_j f,f\rangle|.
\end{align*}
Therefore, $\{A_j\}_{j\in J}$ is $p$-Bessel if and only if so is $\{|A_j|\}_{j\in J}$. Furthermore, it is easy to see that for every $j$ and $f,h\in H$ we also have 
\begin{align*}
    |\beta_j^{|A_j|}(f)-\beta_j^{|A_j|}(h)| = |\beta_j^{A_j}(f)-\beta_j^{A_j}(h)|,
\end{align*}
and 
\begin{align}
\label{eq:passing_to_moduli}
    |\alpha_j^{|A_j|}(f)-\alpha_j^{|A_j|}(h)| = |\alpha_j^{A_j}(f)-\alpha_j^{A_j}(h)|,
\end{align}
Since the moduli of differences of measurements are not changed between $\{A_j\}_{j\in J}$ and $\{|A_j|\}_{j\in J}$, we obtain the following.

\begin{lemma}
\label{lem:negative_iff_positive_failure}
    Let $p\in [2,\infty)$ and let $\{A_j\}_{j\in J}\subseteq B(H)$ be a $p$-Bessel system of semidefinite operators. Then $\{A_j\}_{j\in J}$ does $(p,\alpha)$-SPR with constant $C>0$ if and only if $\{|A_j|\}_{j\in J}$ does $(p,\alpha)$-SPR with constant $C$. The same is true for the map $\beta$.
\end{lemma}

We are ready to prove Theorem \ref{thm:alpha_beta_semidefinite_failure}.

\begin{proof}[Proof of Theorem \ref{thm:alpha_beta_semidefinite_failure}]
    Let $H$ be infinite-dimensional and let $\{A_j\}_{j\in J}\subseteq B(H)$ be a $4$-Bessel system consisting of semidefinite operators. From Proposition \ref{prop:beta_implies_alpha} it suffices to show that $(4,\alpha)$-SPR is impossible. By Lemma \ref{lem:negative_iff_positive_failure}, we may assume $A_j\geq 0$ for all $j$ and by Lemma \ref{lem:finite_systems_in_infinite_Hilbert} we may assume that $|J|=\infty$.

    Fix a unit vector $u\in H$. Since
        $$\sum_{j\in J} \langle A_ju,u\rangle^{2}<\infty,$$
    for every $\varepsilon>0$ we may choose $m$ such that
        $$\sum_{j>m}\langle A_ju,u\rangle^{2}<\varepsilon.$$
    Consider the finite-dimensional subspace
        $$M:=\operatorname{span}\{u,A_1u,\dots,A_mu\}.$$
    Since $H$ is infinite-dimensional, there exists a unit vector $v\in M^\perp$. Thus
        $$v\perp u,\qquad \langle A_ju,v\rangle=0, \quad j=1,\dots,m.$$
    Set
        $$f:=u+v, \qquad h:=u-v.$$
    We will show that for sufficiently small $\varepsilon$, the failure of $(4,\alpha)$-SPR occurs on $f$ and $h$. First, since $u$ and $v$ are orthogonal unit vectors,
        $$\min_{|\theta|=1}\|f-\theta h\|=2.$$
    Next, fix $j\leq m$. Since $\langle A_ju,v\rangle=0$, and since $A_j$ is self-adjoint, we get $\langle A_jv,u\rangle=0$. Using this gives
    \begin{align*}
        \alpha_j(f)^2&=\langle A_j(u+v),u+v\rangle\\
        &=\langle A_ju,u\rangle+\langle A_jv,v\rangle\\
        &=\langle A_j(u-v),u-v\rangle=\alpha_j(h)^2.
    \end{align*}
    Hence
        $$\alpha_j(f)=\alpha_j(h),\qquad j\leq m.$$
    For $j>m$, put $T_j=A_j^{1/2}$. Then $\alpha_j(x)=\|T_jx\|,$ and the reverse triangle inequality gives
    \begin{align*}
        |\alpha_j(f)-\alpha_j(h)|&=\Big|\|T_ju+T_jv\|-\|T_ju-T_jv\|\Big|\\
        &\leq 2\|T_ju\|.
    \end{align*}
    Therefore
    \begin{align*}
    \|\alpha(f)-\alpha(h)\|_4^4&=\sum_{j>m} |\alpha_j(f)-\alpha_j(h)|^4\leq 2^4\sum_{j>m}\|T_ju\|^4\\
    &=2^4\sum_{j>m}\langle A_ju,u\rangle^{2}<2^4\varepsilon.
    \end{align*}
    Consequently,
        $$\|\alpha(f)-\alpha(h)\|_4<2\varepsilon^{1/4}.$$
    Given $C>0$, choose $0<\varepsilon<C^{-4}$. Then 
        $$\min_{|\theta|=1}\|f-\theta h\|=2> C\|\alpha(f)-\alpha(h)\|_4,$$
    showing failure of $(p,\alpha)$-SPR.
\end{proof}

\begin{remark}
    Reproducing the proof of Theorem \ref{thm:alpha_beta_semidefinite_failure} with the exponent $p$ instead of $4$ gives a direct argument why no semidefinite system on an infinite-dimensional Hilbert space can do $(p,\alpha)$-SPR or $(p,\beta)$-SPR.
\end{remark}

Since every self-adjoint rank-one operator is semidefinite, from Theorem \ref{thm:alpha_beta_semidefinite_failure} we can deduce that no $4$-Bessel system of self-adjoint rank-one operators can do $(4,\alpha)$-SPR. Our next goal is to show that $(4,\beta)$-SPR, and thus $(4,\alpha)$-SPR can be achieved on a system of rank-one operators, provided not all of them are self-adjoint. For this purpose, we recall the following inequality. For $f,h\in H$ we have 
\begin{equation}
\label{eq:min_max_Hilbert_Schmidt}
    \min_{|\theta|=1}\|f-\theta h\|\max_{|\theta|=1}\|f-\theta h\|\leq \sqrt{2}\|ff^*-hh^*\|_{\mathrm{HS}},
\end{equation}
where $\|\cdot\|_{\mathrm{HS}}$ denotes the Hilbert-Schmidt norm of an operator. Indeed, note that
    \begin{align*}
        \min_{|\theta|=1}\|f-\theta h\|^2 &= \|f\|^2+\|h\|^2-2|\langle f,h\rangle|, \\
        \max_{|\theta|=1}\|f-\theta h\|^2 &= \|f\|^2+\|h\|^2+2|\langle f,h\rangle|.
    \end{align*}
Moreover, it is standard that 
    $$\|ff^*-hh^*\|_{\mathrm{HS}}^2 = \|f\|^4+\|h\|^4-2|\langle f,h\rangle|^2.$$
Therefore
    $$2\|ff^*-hh^*\|_{\mathrm{HS}}^2-\min_{|\theta|=1}\|f-\theta h\|^2\max_{|\theta|=1}\|f-\theta h\|^2=(\|f\|^2-\|h\|^2)^2\geq 0.$$
Taking the square root gives \eqref{eq:min_max_Hilbert_Schmidt}.

\begin{proof}[Proof of Proposition \ref{prop:beta_alpha_stability_on_rank_1}]
    Let $H$ be infinite-dimensional. We show the existence of a $4$-Bessel system with Bessel bound $B=1$ doing $(4,\beta)$-SPR with constant $\sqrt{2}$. By Proposition \ref{prop:beta_implies_alpha}, the same system does $(4,\alpha)$-SPR with constant $2$.

    Let $(e_n)_{n=1}^\infty\subseteq H$ be an ONB. For $m,n\in\N$ define $E_{mn}\in B(H)$ via
        $$E_{mn}x:=\langle x,e_n\rangle e_m.$$
    Every $E_{mn}$ has rank-one. When $n\neq m$, $E_{mn}^*=E_{nm} \neq E_{mn}$. We show now that the system $\{E_{mn}\}_{m,n\in \N}$ is $4$-Bessel. For every $f\in H$ we have
    \begin{align*}
        \sum_{m,n}|\langle E_{mn}f,f\rangle|^2 &= \sum_{m,n}|\langle f,e_m\rangle|^2 \cdot |\langle f, e_n\rangle|^2 \\
        &= \left(\sum_m |\langle f,e_m\rangle|^2\right)\cdot \left(\sum_n |\langle f,e_n\rangle|^2\right) \\
        &= \|f\|^2\cdot \|f\|^2 = \|f\|^4.
    \end{align*}
    Next, observe that for every $m,n\in \N$ and $f\in H$ we have
        $$\beta_{mn}(f)=\langle f,e_n\rangle\cdot\overline{\langle f,e_m\rangle}.$$
    Using the definition that $(ff^*)x=\langle x,f\rangle f$ for every $x\in H$, we get
        $$(ff^*)e_m=\langle e_m,f\rangle f = \overline{\langle f,e_m\rangle}f.$$
    Therefore
        $$\langle (ff^*)e_m,e_n\rangle = \overline{\langle f,e_m\rangle}\langle f,e_n\rangle=\beta_{mn}(f)$$
    Hence 
    \begin{align*}
        \beta_{mn}(f)-\beta_{mn}(h)&=\langle (ff^*)e_m,e_n\rangle - \langle (hh^*)e_m,e_n\rangle \\
        &= \langle (ff^*-hh^*)e_m,e_n\rangle.
    \end{align*}
    Using the fact that the Hilbert-Schmidt norm is the $\ell_2$-norm of all operator entries, we get 
    \begin{align}
    \label{eq:beta_is_Hilbert_Schmidt}
    \|\beta(f)-\beta(h)\|_2^2 &= \sum_{m,n} |\langle (ff^*-hh^*)e_m,e_n\rangle |^2=\|ff^*-hh^*\|_{\mathrm{HS}}^2.
    \end{align}
    Hence, by \eqref{eq:min_max_Hilbert_Schmidt} we obtain
        $$\min_{|\theta|=1}\|f-\theta h\|\max_{|\theta|=1}\|f-\theta h\|\leq \sqrt{2}\|ff^*-hh^*\|_{\mathrm{HS}}=\sqrt{2}\|\beta(f)-\beta(h)\|_2.$$
    Therefore $\{E_{mn}\}_{m,n\in \N}$ does $(4,\beta)$-SPR with $C=\sqrt{2}$.
\end{proof}

\begin{example}
    When $H=\C^2$, then given the standard ONB $\{e_1,e_2\}$, for every $f=ae_1+be_2\in H$ we get 
    $$\beta(f)=(E_{11}(f),E_{12}(f),E_{21}(f), E_{22}(f))=(|a|^2,a\overline{b}, \overline{a}b,|b|^2).$$
    These measurements uniquely determine $f$, up to a unimodular scalar.
\end{example}

\begin{remark}
     We note that the system from Proposition \ref{prop:beta_alpha_stability_on_rank_1} is $p$-Bessel if and only if $p\in [4,\infty)$. Indeed, since it is $4$-Bessel, according to Lemma \ref{lem:p_Bessel_relation} it is $p$-Bessel for $p\in [4,\infty)$. On the other hand, for $p\in [2,4)$, consider 
        $$f_N=\frac{1}{\sqrt{N}}\sum_{j=1}^N e_j, \qquad N\geq 1.$$
    Then
        $$\sum_{m,n}|\langle E_{mn}f_N,f_N\rangle|^{p/2}=N^2N^{-p/2}=N^{2-p/2}\xrightarrow{N\to\infty}\infty.$$
    Therefore, this system is not $p$-Bessel for $p\in [2,4)$. In particular, the measurement maps $\beta$ and $\alpha$ are not well-defined in those ranges. From Theorem \ref{thm:alpha_beta_failure_on_strong_Bessel} we also know that for $p>4$, $\{E_{mn}\}_{m,n\in \N}$ fails $(p,\alpha)$-SPR. This can be observed directly, by considering the pairs $(f_N,0)$ for $N\geq 1$. We have that 
        $$\min_{|\theta|=1}\|f_N-\theta\cdot 0\|= \|f_N\|=1, \qquad N\geq 1,$$
    while
        $$\|\alpha(f_N)-\alpha(0)\|_p = \left(N^2N^{-p/2}\right)^{1/p}=N^{2/p-1/2}\xrightarrow{N\to\infty}0.$$
    Therefore, $\{E_{mn}\}_{m,n\in \N}\subseteq B(H)$ does $(p,\alpha)$-SPR if and only if $p=4$. 
\end{remark}

To achieve $\beta$-stability (and therefore also $\alpha$-stability) on an indefinite system of self-adjoint operators, one has to assume that it does not consist entirely of operators of rank at most one. 

\begin{proof}[Proof of Theorem \ref{thm:4_Bessel_stability}]
    Let $H$ be infinite-dimensional. First, we show the existence of a self-adjoint $4$-Bessel system of operators of rank at most two doing $(4,\beta)$-SPR with constant $\sqrt{2}$.
    
    Let $(e_n)_{n\in \N}$ be an ONB of $H$. As in the proof of Proposition \ref{prop:beta_alpha_stability_on_rank_1}, define
        $$E_{mn}x:=\langle x,e_n\rangle e_m, \qquad x\in H, \quad n,m\in \N.$$
    Next, consider the self-adjoint operators given by
        $$S_{mn}:=\frac{1}{\sqrt{2}}(E_{mn}+E_{nm}), \qquad m<n.$$
    If $H$ is complex, we additionally take
        $$T_{mn}:=\frac{i}{\sqrt{2}}(E_{mn} - E_{nm}), \qquad m<n.$$
    In the real case, the family
        $$\{E_{nn}:n\geq1\}\cup\{S_{mn}:m<n\}$$
    is an ONB of the real Hilbert space of Hilbert-Schmidt self-adjoint operators
    $\mathcal{S}_2(H)_{\mathrm{sa}}$. In the complex case, the family
        $$\{E_{nn}:n\geq1\}\cup\{S_{mn},T_{mn}:m<n\}$$
    is an ONB of $\mathcal{S}_2(H)_{\mathrm{sa}}$, regarded as a real Hilbert space. Enumerate the corresponding orthonormal basis as $\{A_j\}_{j\in J}$. Each $A_j$ is self-adjoint and has rank either $1$ or $2$. By Parseval's identity
        $$\sum_{j\in J} |\langle A_jf,f\rangle|^2=\|ff^*\|_{\mathrm{HS}}^2=\|f\|^4.$$
    Hence $\{A_j\}_{j\in J}$ is $4$-Bessel with Bessel bound $1$. Moreover, by \eqref{eq:beta_is_Hilbert_Schmidt}
        $$\|\beta(f)-\beta(h)\|_2=\|ff^*-hh^*\|_{\mathrm{HS}}.$$
    Therefore,  \eqref{eq:min_max_Hilbert_Schmidt} implies
        $$\min_{|\theta|=1}\|f-\theta h\|\max_{|\theta|=1}\|f-\theta h\|\leq\sqrt{2}\|\beta(f)-\beta(h)\|_{2},$$
    so the system $\{A_j\}_{j\in J}$ does $(4,\beta)$-SPR with constant $\sqrt{2}$. By Proposition \ref{prop:beta_implies_alpha}, it also does $(4,\alpha)$-SPR with constant $2$.

    Next, we fix an integer $r\geq 2$. Our goal is to perturb $\{A_j\}_{j\in J}$ to construct a $4$-Bessel system in which every operator has rank exactly $r$. Without loss of generality we may assume that $J=\N$. First, note that for each $j\in \N$ there exists a finite set $S_j\subseteq \N$, such that 
        $$\operatorname{Range}(A_j)\subseteq \lspan\{e_k: k\in S_j\}.$$ 
    Next, pick any set $F_1\subseteq \N\setminus S_1$ satisfying $|F_1|=r-\operatorname{rank}(A_1).$ We note that 
        $$\lspan\{e_k: k\in F_1\} \perp \operatorname{Range}(A_1).$$
    Recursively, we pick sets $F_j$ for $j\geq 2$ satisfying
        $$F_j\subseteq\mathbb N\setminus\left(S_j\cup\bigcup_{i<j}F_i\right), \qquad |F_j|=r-\operatorname{rank}(A_j).$$
    By construction, $F_j$ are pairwise disjoint and moreover 
        $$\operatorname{span}\{e_k:k\in F_j\}\perp\operatorname{Range}(A_j), \qquad j\in \N.$$
    If $r-\operatorname{rank}(A_j)>0$, define
        $$D_j:=\frac{1}{\sqrt{|F_j|}}\sum_{k\in F_j}E_{kk},$$
    and if $r-\operatorname{rank}(A_j)=0$, set $D_j:=0$. Then every $D_j$ is self-adjoint,
        $$\operatorname{rank}(D_j)=|F_j|=r-\operatorname{rank}(A_j),$$
    and
        $$\operatorname{Range}(D_j)\perp\operatorname{Range}(A_j).$$
    Moreover, the nonzero operators $D_j$ form an orthonormal sequence in $\mathcal{S}_2(H)_{\mathrm{sa}}$.
    
    Fix $0<\varepsilon<1$ and for each $j\in \N$ define
        $$B_j:=A_j+\varepsilon D_j.$$
    We claim that the system $\{B_j\}_{j\in \N}$ satisfies the statement of the theorem. We need to show the following three claims:
    \begin{enumerate}
        \item for each $j\in \N$, $\operatorname{rank}(B_j)=r;$
        \item $\{B_j\}_{j\in \N}$ is $4$-Bessel with Bessel bound $(1+\varepsilon)^2$;
        \item $\{B_j\}_{j\in\N}$ does $\beta$-SPR with constant $\frac{\sqrt{2}}{1-\varepsilon}$.
    \end{enumerate}    
    First, since $\operatorname{Range}(A_j)$ and $\operatorname{Range}(D_j)$ are orthogonal, self-adjointness and orthogonality of the ranges imply $A_jD_j=D_jA_j=0$. Hence
        $$\operatorname{rank}(B_j)=\operatorname{rank}(A_j)+\operatorname{rank}(D_j)=r,$$
    which gives (1).
    
    We now show (2). Fix $f\in H$. By Parseval's identity we have
        $$\left(\sum_{j\in\N}|\langle ff^*,A_j\rangle_{\mathrm{HS}}|^2\right)^{1/2} = \|ff^*\|_{\mathrm{HS}}.$$
    Applying Bessel's inequality we also obtain
        $$\left(\sum_{j\in\N}|\langle ff^*,D_j\rangle_{\mathrm{HS}}|^2\right)^{1/2} \leq \|ff^*\|_{\mathrm{HS}}.$$
    Combining these with the triangle inequality gives
    \begin{align*}
    \left(\sum_{j\in\N}|\langle B_jf,f\rangle|^2\right)^{1/2}&=\left(\sum_{j\in\N}|\langle ff^*,B_j\rangle_{\mathrm{HS}}|^2\right)^{1/2}\\
    &\leq\left(\sum_{j\in\N}|\langle ff^*,A_j\rangle_{\mathrm{HS}}|^2\right)^{1/2}+\varepsilon\left(\sum_{j\in\N}|\langle ff^*,D_j\rangle_{\mathrm{HS}}|^2\right)^{1/2} \\
    &\leq (1+\varepsilon)\|ff^*\|_{\mathrm{HS}}=(1+\varepsilon)\|f\|^2.
    \end{align*}
    By squaring both sides we obtain 
        $$\sum_{j\in\N}|\langle B_jf,f\rangle|^2\leq (1+\varepsilon)^2\|f\|^4,$$
    which establishes (2). Finally, fix $f,h\in H$. By repeating similar calculations as above for $ff^*-hh^*$ and using the reverse triangle inequality we have
    \begin{align*}
        \|\beta(f)-\beta(h)\|_2&=\left(\sum_{j\in\N} |\langle ff^*-hh^*, B_j\rangle_{\mathrm{HS}}
        |^2\right)^{1/2}\\
        &\geq \left(\sum_{j\in\N}|\langle ff^*-hh^*,A_j\rangle_{\mathrm{HS}}|^2\right)^{1/2}-\varepsilon\left(\sum_{j\in\N}|\langle ff^*-hh^*,D_j\rangle_{\mathrm{HS}}|^2\right)^{1/2} \\
        &\geq (1-\varepsilon)\|ff^*-hh^*\|_{\mathrm{HS}}.
    \end{align*}
    Using \eqref{eq:min_max_Hilbert_Schmidt}, we therefore get 
        $$\min_{|\theta|=1}\|f-\theta h\|\max_{|\theta|=1}\|f-\theta h\|\leq \frac{\sqrt{2}}{1-\varepsilon}\|\beta(f)-\beta(h)\|_2.$$
    Hence $\{B_j\}_{j\in\N}$ does $(4,\beta)$-SPR with constant $C=\frac{\sqrt{2}}{1-\varepsilon}.$ Applying Proposition \ref{prop:beta_implies_alpha} gives the second part of the theorem for the map $\alpha$.
    \end{proof}

\section{Orthogonal reductions and finite-dimensional stability}
\label{sec:orthogonal_reductions_and_f_dim_stability}

In this section $H$ stands for a finite or infinite-dimensional Hilbert space, unless stated otherwise. Recall that a frame $\{g_j\}_{j\in J}\subseteq H$ does SPR with constant $C>0$ if for every pair $f,h\in H$, we have 
$$\min_{|\theta|=1}\|f-\theta h\|\leq C\|\mathcal{A}(f)-\mathcal{A}(h)\|_2.$$
It is a known result in classical phase retrieval that in order to check if the stability inequality is satisfied for a frame $\{g_j\}_{j\in J}$, it suffices to verify stability inequality on orthogonal pairs. More precisely, 
\cite[Theorem 2.1]{Alharbi2024} shows that $\{g_j\}_{j\in J}$ does SPR with constant $C>0$ if and only if for every orthogonal pair $f,h\in H$ we have
$$\min_{|\theta|=1}\|f-\theta h\|\leq C\|\mathcal{A}(f)-\mathcal{A}(h)\|_2.$$

A related uniqueness result was obtained in the setting of orthogonal projections in $\R^d$ \cite[Lemma 4.1]{Cahill2016}; a finite system of orthogonal projections $\{A_j\}_{j\in J}\subseteq B(H)_+$ does $\alpha$-PR if and only if for every orthogonal $f,h\in H$ we have 
$$\alpha(f)=\alpha(h) \implies f=\theta h, \text{ for some }|\theta|=1.$$

First, we will show that such a reduction is always possible for the maps $\beta$ and $\alpha$. We will also show that the stability constant for $\beta$ is always preserved, whereas for $\alpha$ it is preserved whenever we are given a system of semidefinite operators. 

We begin with the following auxiliary lemma which will be an important tool in the orthogonal reduction performed later.
\begin{lemma}
\label{lem:orthogonal_reduction_2}
    Let $x,y\in H$. Then for every $t\in [0,1]$ we have 
    $$\Big|\|tx+y\| - \|tx-y\| \Big|\leq \Big|\|x+y\| - \|x-y\| \Big|.$$
\end{lemma}
\begin{proof}
    Fix $x,y\in H$. If $x=0$ or $y=0$, the inequality follows immediately, so assume $x,y\neq 0$. Define $\varphi:[0,1]\to [0,\infty)$ by
        $$\varphi(t):=\Big|\|tx+y\| - \|tx-y\| \Big|.$$
    Our goal is to show that $\varphi$ is non-decreasing on $[0,1]$, which will finish the proof. Since $\varphi$ is non-negative, it suffices to show that $\varphi^2$ is non-decreasing on $[0,1]$.

    For this purpose, denote $a=\|x\|^2, b=\|y\|^2$ and $c=\re\langle x,y\rangle$. Note that by the Cauchy-Schwarz inequality we have $ab-c^2\geq 0$. A direct calculation gives that 
        $$\varphi^2(t)=2(at^2+b) - 2\sqrt{(at^2+b)^2 - 4t^2c^2}.$$
    Note that for every $t\in [0,1]$, from $c^2\leq ab$ we have 
        $$(at^2+b)^2-4t^2c^2 \geq (at^2+b)^2 - 4abt^2 = (at^2-b)^2\geq 0.$$
    Fix $t\in [0,1]$ and suppose that $D(t):=(at^2+b)^2-4t^2c^2 >0$. Then 
        $$\frac{d}{dt}\varphi^2(t)=4t\Big(a-\frac{a(at^2+b)-2c^2}{\sqrt{D(t)}}\Big).$$
    If $a(at^2+b)-2c^2\leq 0$, then $\frac{d}{dt}\varphi^2(t)\geq 0$. Suppose now that $a(at^2+b)-2c^2> 0$. A direct calculation shows that 
        $$a^2D(t)-[a(at^2+b)-2c^2]^2 = 4c^2(ab-c^2)\geq 0.$$
    As both $D(t)$ and $a(at^2+b)-2c^2$ are positive, taking the square root gives 
        $$a\geq \frac{a(at^2+b)-2c^2}{\sqrt{D(t)}}.$$
    Therefore again $\frac{d}{dt}\varphi^2(t)\geq 0$. 

    This shows that $\varphi^2(t)$ is non-decreasing on every interval where $D(t)>0$. When $D(t)=0$, since $D(t)=\|tx+y\|^2\|tx-y\|^2$, hence $D$ vanishes at most once on $[0,1]$. Therefore continuity of $\varphi^2$ forces it to be non-decreasing everywhere. This finishes the proof.
\end{proof}

The following lemma will provide a way of reducing the study of $(p,\alpha)$-SPR and $(p,\beta)$-SPR to the study of orthogonal pairs.

\begin{lemma}
\label{lem:orthogonal_reduction}
    Let $\{A_j\}_{j\in J}\subseteq B(H)$ and let $f,h\in H$ be linearly independent. There exist non-zero orthogonal $f_0,h_0\in \lspan\{f,h\}$ and $t\in (0,1]$ such that
    \begin{enumerate}
        \item $\min_{|\theta|=1}\|f_0-\theta h_0\|\max_{|\theta|=1}\|f_0-\theta h_0\|=t\cdot \min_{|\theta|=1}\|f-\theta h\|\max_{|\theta|=1}\|f-\theta h\|;$
        \item $\min_{|\theta|=1}\|f_0-\theta h_0\|= \min_{|\theta|=1}\|f-\theta h\|$;
        \item for every $j$ 
        $$\beta_j(f_0)-\beta_j(h_0)=t\cdot (\beta_j(f)-\beta_j(h)).$$
    \end{enumerate}
    If in addition $\{A_j\}_{j\in J}\subseteq B(H)$ is a system of semidefinite operators, for every $j$ we also have
        $$|\alpha_j(f_0)-\alpha_j(h_0)|\leq |\alpha_j(f)-\alpha_j(h)|.$$
\end{lemma}
\begin{proof}
    Fix linearly independent $f,h\in H$. Since replacing $h$ by a unimodular multiple does not change the $\alpha$ or $\beta$ measurements, we may assume that $\langle f,h\rangle\geq 0$. Set 
    $$u=\frac{f+h}{2},\qquad v=\frac{f-h}{2}.$$
    Then $f=u+v$, $h=u-v$ and moreover $\|u\|^2 - \|v\|^2 = \langle f,h\rangle \geq 0$. We also note that 
        $$\langle u,v\rangle = \frac{1}{4}\langle f+h,f-h\rangle = \frac{1}{4}(\|f\|^2-\|h\|^2)\in \R.$$
    Next, linear independence of $f$ and $h$ guarantees that $u$ and $v$ are also linearly independent, in particular $u,v\neq 0$. Therefore $t:=\frac{\|v\|}{\|u\|}\in (0,1]$ is well-defined. Set 
        $$f_0 := tu+v, \qquad h_0:=tu-v.$$
    As $\langle u,v\rangle\in\R$, we get $\langle f_0,h_0\rangle = t^2\|u\|^2 - \|v\|^2 = 0$, thus $f_0\perp h_0$. Moreover, linear independence of $u$ and $v$ guarantees that $f_0,h_0\neq 0$. Next, for each $j$,
        $$\beta_j(f_0)-\beta_j(h_0) = 2t(\langle A_ju,v\rangle + \langle A_jv,u\rangle) = t\cdot (\beta_j(f)-\beta_j(h)).$$
    We also have
    $$\min_{|\theta|=1}\|f-\theta h\|\max_{|\theta|=1}\|f-\theta h\| = 4\|u\|\|v\|,$$ 
    whereas
        $$\min_{|\theta|=1}\|f_0-\theta h_0\|\max_{|\theta|=1}\|f_0-\theta h_0\| = \|f_0\|^2+\|h_0\|^2 = 4\|v\|^2 = 4t\|u\|\|v\|.$$
    Hence 
        $$\min_{|\theta|=1}\|f_0-\theta h_0\|\max_{|\theta|=1}\|f_0-\theta h_0\|=t\cdot \min_{|\theta|=1}\|f-\theta h\|\max_{|\theta|=1}\|f-\theta h\|.$$
    Furthermore, we also get
        $$\min_{|\theta|=1}\|f_0-\theta h_0\|=2\|v\|=\min_{|\theta|=1}\|f-\theta h\|.$$
    Suppose now that $A_j$ is semidefinite for each $j$. By \eqref{eq:passing_to_moduli}, we can assume that $A_j\geq 0$ for all $j$. Fix $j$ and set $T_j=A^{1/2}_j$. Then $\alpha_j(f)=\|T_jf\|$. Applying Lemma \ref{lem:orthogonal_reduction_2} with $x=T_ju$ and $y=T_jv$ gives
    \begin{align*}
        |\alpha_j(f_0)-\alpha_j(h_0)| &=\Big|\|tT_ju + T_jv\| - \|tT_ju-T_jv\| \Big|\\
        &\leq \Big|\|T_ju + T_jv\| - \|T_ju - T_jv\|\Big|\\
        &=|\alpha_j(f) - \alpha_j(h)|.
    \end{align*}
\end{proof}

We are ready to show that for the map $\beta$, orthogonal reduction is not only possible, but always preserves the stability constant.

\begin{proof}[Proof of Theorem \ref{thm:beta_sufficient_on_ortho_pairs}]
    Let $\{A_j\}_{j\in J}\subseteq B(H)$ be a $p$-Bessel system for some $p\in [2,\infty)$. If $\{A_j\}_{j\in J}$ does $\beta$-SPR with constant $C>0$, the $\beta$-SPR inequality is also satisfied on orthogonal pairs. To prove the converse, suppose that for every orthogonal pair $f_0,h_0\in H$ we have
    $$\min_{|\theta|=1}\|f_0-\theta h_0\|\max_{|\theta|=1}\|f_0-\theta h_0\|\leq C\|\beta(f_0)-\beta(h_0)\|_{p/2}.$$
    Our goal is to show that the same holds for any pair $f,h\in H$. Fix such a pair. 

    \textbf{Case 1.} $f=ch$ for some $c\in \F$. Then 
        $$\beta(f)-\beta(h)=|c|^2\beta(h) - \beta(h)=(|c|^2-1)\beta(h)$$ 
    and 
        $$\min_{|\theta|=1}\|f-\theta h\|\max_{|\theta|=1}\|f-\theta h\| = ||c|^2-1|\|h\|^2.$$
    Applying the $(p,\beta)$-SPR inequality to the orthogonal pair $(0,h)$ gives $\|h\|^2\leq C\|\beta(h)\|_{p/2}$. Therefore 
    \begin{align*}
        \min_{|\theta|=1}\|f-\theta h\|\max_{|\theta|=1}\|f-\theta h\| &= |1-|c|^2|\|h\|^2\leq C|1-|c|^2|\|\beta(h)\|_{p/2} \\
        &=C\|\beta(f)-\beta(h)\|_{p/2}.
    \end{align*}
    A similar argument can be applied to the case $h=0$ and $f\neq 0$.
    
    \textbf{Case 2.} Suppose that $f$ and $h$ are linearly independent. From Lemma \ref{lem:orthogonal_reduction} find orthogonal $f_0,h_0\in \lspan\{f,h\}$ and $t\in (0,1]$ such that 
    $$\min_{|\theta|=1}\|f_0-\theta h_0\|\max_{|\theta|=1}\|f_0-\theta h_0\|=t\cdot \min_{|\theta|=1}\|f-\theta h\|\max_{|\theta|=1}\|f-\theta h\|$$
    and 
    $$\beta(f_0)-\beta(h_0)=t\cdot (\beta(f)-\beta(h)).$$
    We get
    \begin{align*}
        t\cdot \min_{|\theta|=1}\|f-\theta h\|\max_{|\theta|=1}\|f-\theta h\|&=\min_{|\theta|=1}\|f_0-\theta h_0\|\max_{|\theta|=1}\|f_0-\theta h_0\|\\
        &\leq C\|\beta(f_0)-\beta(h_0)\|_{p/2} \\
        &=t\cdot C\|\beta(f)-\beta(h)\|_{p/2}.
    \end{align*}
    Since $t\neq 0$, we get 
    $$\min_{|\theta|=1}\|f-\theta h\|\max_{|\theta|=1}\|f-\theta h\|\leq C\|\beta(f)-\beta(h)\|_{p/2}.$$
\end{proof}

Now we prove a similar orthogonality reduction for the map $\alpha$. For systems of semidefinite operators, we show that this reduction also preserves the stability constant.

\begin{proof}[Proof of Theorem \ref{thm:alpha_sufficient_on_ortho_pairs}]
    Let $\{A_j\}_{j\in J}\subseteq B(H)$ be a $p$-Bessel system with constant $B$ for some $p\in [2,\infty)$. If $\{A_j\}_{j\in J}$ does $\alpha$-SPR with constant $C>0$, the $\alpha$-SPR inequality is also satisfied on orthogonal pairs. To prove the converse, suppose that for every orthogonal pair $f_0,h_0\in H$ we have
    $$\min_{|\theta|=1}\|f_0-\theta h_0\|\leq C\|\alpha(f_0)-\alpha(h_0)\|_p.$$
    We claim that $\{A_j\}_{j\in J}$ does $(p,\beta)$-SPR on orthogonal pairs with constant $2C^2$. Indeed, fix orthogonal $f_0,h_0\in H$. Since $\max_{|\theta|=1}\|f_0-\theta h_0\|=\min_{|\theta|=1}\|f_0-\theta h_0\|$, we have 
    $$\max_{|\theta|=1}\|f_0-\theta h_0\|\min_{|\theta|=1}\|f_0-\theta h_0\|\leq C^2\|\alpha(f_0)-\alpha(h_0)\|_p^2.$$
    Moreover, using Lemma \ref{lem:inequality_for_alpha} we get 
        $$\|\alpha(f_0)-\alpha(h_0)\|_p^2\leq 2\|\beta(f_0)-\beta(h_0)\|_{p/2}.$$
    Combining these inequalities gives 
    \begin{align*}
        \max_{|\theta|=1}\|f_0-\theta h_0\|\min_{|\theta|=1}\|f_0-\theta h_0\|&\leq C^2\|\alpha(f_0)-\alpha(h_0)\|_p^2\\
        &\leq 2C^2\|\beta(f_0)-\beta(h_0)\|_{p/2},
    \end{align*}
    which proves the claim. From Theorem \ref{thm:beta_sufficient_on_ortho_pairs}, the system $\{A_j\}_{j\in J}$ does $(p,\beta)$-SPR with constant $2C^2$. From Proposition \ref{prop:beta_implies_alpha}, this implies that $\{A_j\}_{j\in J}$ does $(p,\alpha)$-SPR with constant $2\sqrt{2}C^2B^{1/p}$.

    Suppose now that in addition $\{A_j\}_{j\in J} \subseteq B(H)$ is a system of semidefinite operators. Since again the forward direction is immediate, suppose that for every orthogonal $f_0,h_0\in H$ we have 
        $$\min_{|\theta|=1}\|f_0-\theta h_0\|\leq C\|\alpha(f_0)-\alpha(h_0)\|_p.$$
    We will show that the same inequality holds for any pair $f,h\in H$. Fix such a pair. If $f=ch$ for some $c\in \F$, the $(p,\alpha)$-SPR inequality with constant $C$ can be deduced via the orthogonal pair $(0,h)$, as in the proof of Theorem \ref{thm:beta_sufficient_on_ortho_pairs}. Therefore suppose that $f$ and $h$ are linearly independent. From Lemma \ref{lem:orthogonal_reduction} we can find non-zero orthogonal $f_0,h_0\in H$ such that 
        $$\min_{|\theta|=1}\|f-\theta h\| = \min_{|\theta|=1}\|f_0-\theta h_0\|$$
    and
        $$\|\alpha(f_0) - \alpha(h_0)\|_p\leq \|\alpha(f)-\alpha(h)\|_p.$$
    Thus
    \begin{align*}
        \min_{|\theta|=1}\|f-\theta h\| &= \min_{|\theta|=1}\|f_0-\theta h_0\| \leq C\|\alpha(f_0) - \alpha(h_0)\|_p \\
        &\leq C\|\alpha(f)-\alpha(h)\|_p.
    \end{align*}
    Hence $\{A_j\}_{j\in J}$ does $(p,\alpha)$-SPR with constant $C$.
\end{proof}

As we will see in Example \ref{ex:stability_constant_not_achieved_on_ortho}, part (1) of Theorem \ref{thm:alpha_sufficient_on_ortho_pairs} cannot be improved in general, as the optimal global stability constant may be strictly larger than the optimal constant on orthogonal pairs.

The first important consequence of the orthogonal reduction allows us to show that $(p,\alpha)$-SPR and $(p,\beta)$-SPR are equivalent up to a constant.

\begin{proof}[Proof of Corollary \ref{cor:beta_and_alpha_are_equivalent}]
    Let $\{A_j\}_{j\in J}\subseteq B(H)$ be $p$-Bessel. If $\{A_j\}_{j\in J}$ does $(p,\beta)$-SPR, then it does $(p,\alpha)$-SPR by Proposition \ref{prop:beta_implies_alpha}.

    Suppose now that $\{A_j\}_{j\in J}$ does $(p,\alpha)$-SPR with constant $C>0$. Our goal is to show it does $(p,\beta)$-SPR. By Theorem \ref{thm:beta_sufficient_on_ortho_pairs} it suffices to check the $(p,\beta)$-SPR inequality on orthogonal pairs. Fix orthogonal $f,h\in H$. Then 
    \begin{align*}
    \min_{|\theta|=1}\|f-\theta h\|\max_{|\theta|=1}\|f-\theta h\| &= \min_{|\theta|=1}\|f-\theta h\|^2 \leq C^2\|\alpha(f)-\alpha(h)\|_p^2 \\
    &\leq 2C^2\|\beta(f)-\beta(h)\|_{p/2}.
    \end{align*}
    Hence $\{A_j\}_{j\in J}$ does $(p,\beta)$-SPR.
\end{proof}

A second consequence of the orthogonal reduction is that in finite-dimensional setting, stability for $\alpha$ and $\beta$ maps is achieved whenever they perform phase retrieval, just like in the classical phase retrieval.

\begin{proof}[Proof of Theorem \ref{thm:beta_and_alpha_PR_in_finite_dim}]
    Let $H$ be a finite-dimensional Hilbert space. Since $(p,\alpha)$-SPR and $(p,\beta)$-SPR are equivalent by Corollary \ref{cor:beta_and_alpha_are_equivalent}, and since $\alpha$-PR and $\beta$-PR are equivalent by \eqref{eq:alpha_PR_iff_beta_PR}, it suffices to show that if $\{A_j\}_{j\in J}$ does $\beta$-PR, it does $(p,\beta)$-SPR. By Theorem \ref{thm:beta_sufficient_on_ortho_pairs} it is enough to show that $\beta$ does SPR on orthogonal pairs. 

    Suppose otherwise. Then we can find a pair of sequences $(f_n),(h_n)\subseteq H$ such that $f_n\perp h_n$ for every $n$ and moreover
        $$\|f_n\|^2+\|h_n\|^2=\min_{|\theta|=1}\|f_n-\theta h_n\|\max_{|\theta|=1}\|f_n-\theta h_n\|>n\|\beta(f_n)-\beta(h_n)\|_{p/2}.$$
    Without loss of generality we can assume that for every $n$ we have $\|f_n\|\geq \|h_n\|$. By replacing $f_n$ and $h_n$ with $\frac{f_n}{\|f_n\|}$ and $\frac{h_n}{\|f_n\|}$, we may additionally assume that $\|f_n\|=1\geq \|h_n\|$. This gives
        $$2>n\|\beta(f_n)-\beta(h_n)\|_{p/2}.$$
    Since $(f_n),(h_n)\subseteq B_H$, which is a compact set due to $H$ being finite dimensional, by passing to subsequences if necessary, we can find $f,h\in H$ such that $f_n\to f$ and $h_n\to h$. We claim that $\{A_j\}_{j\in J}$ fails $\beta$-PR on $f,h$. First, observe that for every $j$ we have
        $$|\beta_j(f_n)-\beta_j(h_n)|\leq \|\beta(f_n) - \beta(h_n)\|_{p/2}<\frac{2}{n}.$$
    Since $A_j$ is bounded and $f_n\to f, h_n\to h$, we get $\beta_j(f)=\beta_j(h)$ for every $j$. Therefore $\beta(f)=\beta(h).$
    
    On the other hand, using the fact that $f\perp h$ and $\|f\|=1\geq \|h\|$ we get 
        $$1\leq \|f\|^2+\|h\|^2 = \min_{|\theta|=1}\|f-\theta h\|\max_{|\theta|=1}\|f-\theta h\|$$
    thus $f\neq \theta h$ for all $|\theta|=1$.
\end{proof}

The following gives an example showing that the optimal global constant can be strictly larger than the optimal constant on orthogonal pairs, even if the system consists of self-adjoint operators.

\begin{example}
\label{ex:stability_constant_not_achieved_on_ortho}
    Let $H=\R^2$ and consider the operators
        $$A_1=\begin{pmatrix}
                1&0\\
                0&-1
            \end{pmatrix},
            \qquad
            A_2=\begin{pmatrix}
                0&1\\
                1&0
            \end{pmatrix},
            \qquad
            A_3=I.$$
    Note that each $A_j$ is self-adjoint. Therefore, for every $j=1,2,3$ and $z\in \R^2$ we have $\alpha_j(z)=\sgn(\langle A_jz,z\rangle)\sqrt{|\langle A_jz,z\rangle|}$. 

    We claim the following:
    \begin{enumerate}
        \item $\{A_j\}_{j=1}^3$ does $(2,\alpha)$-SPR with some optimal constant $C_{glob}$;
        \item The stability constant satisfies $C_{glob}>\frac{1}{\sqrt{2}}$;
        \item On orthogonal pairs, $\{A_j\}_{j=1}^3$ satisfies the $(2,\alpha)$-SPR inequality with optimal constant $C_{\perp}=\frac{1}{\sqrt{2}}$.
    \end{enumerate}
    With these three claims, we see that the stability constant is not achieved on orthogonal pairs.
    
    \textbf{Claim (1).} Since $H$ is finite dimensional, from Theorem \ref{thm:beta_and_alpha_PR_in_finite_dim} it is enough to show that $\{A_j\}_{j=1}^3$ does $\alpha$-PR. Fix $z=(x,y)\in\R^2$. A direct calculation gives
        $$\alpha(x,y)=(\sgn(x^2-y^2)\sqrt{|x^2-y^2|}, \sgn(2xy)\sqrt{|2xy|}, \sqrt{x^2+y^2}).$$
    Suppose for $f=(a,b)$ and $h=(c,d)$ in $\R^2$ we have $\alpha(f)=\alpha(h)$. Then $\beta(f)=\beta(h)$. In particular, 
    \begin{align*}
        a^2-b^2 = c^2-d^2, \qquad 2ab = 2cd, \qquad a^2+b^2=c^2+d^2.
    \end{align*}
    Then
    $$a^2 = \frac{(a^2+b^2)+(a^2-b^2)}{2} = \frac{(c^2+d^2)+(c^2-d^2)}{2}=c^2$$
    therefore $|a|=|c|$. Similarly we obtain $|b|=|d|$. Since in addition $ab=cd$, this gives that $(a,b)=\pm (c,d)$. Hence $\{A_j\}_{j=1}^3$ does $\alpha$-PR.

    \textbf{Claim (2).} For $\gamma\in [0,2\pi]$ define $u_\gamma=(\cos\gamma, \sin\gamma)$. Then 
    $$\alpha(u_\gamma) = (\sgn(\cos2\gamma)\sqrt{|\cos 2\gamma|},\sgn(\sin2\gamma)\sqrt{|\sin 2\gamma|},1).$$
    Take $h=u_{\pi/8}$ and $f=u_{\pi/6}$. Then 
        $$\alpha(h) = (2^{-1/4}, 2^{-1/4},1), \qquad \alpha(f)=(2^{-1/2}, 3^{1/4}2^{-1/2}, 1).$$
    Hence 
        $$\|\alpha(f)-\alpha(h)\|_2 = \Big[\big(2^{-1/2} - 2^{-1/4}\big)^2+\big(3^{1/4}2^{-1/2}-2^{-1/4}\big)^2\Big]^{1/2}<0.162.$$
    On the other hand, 
    \begin{align*}
        \min_{|\theta|=1}\|f-\theta h\| &= \sqrt{\|f\|^2+\|h\|^2-2\langle f,h\rangle} = \sqrt{2-2\cos(\pi/24)} \\
        &=2\sin(\pi/48)> 0.13.
    \end{align*}
    Therefore
        $$C_{glob}\geq C_{(f,h)}=\frac{\min_{|\theta|=1}\|f-\theta h\|}{\|\alpha(f)-\alpha(h)\|_2} > \frac{0.13}{0.162} > 0.8> \frac{1}{\sqrt{2}}.$$
    
    \textbf{Claim (3).} Fix any non-zero orthogonal $f,h\in \R^2$. First, write $f=ru_{\gamma}$ for some $\gamma\in[0,2\pi]$ and $r=\|f\|>0$. Since $f\perp h$, $h\in u_{\gamma}^\perp = \lspan\{u_{\gamma+\pi/2}\}$. Thus $h=s u_{\gamma+\pi/2}$ for some $s\neq 0$. Since replacing $h$ by $-h$ does not change the measurement map, we may assume that $s> 0$. Then 
    \begin{align*}
        \alpha(f)-\alpha(h) = ((r+s)\sgn(\cos2\gamma)\sqrt{|\cos2\gamma|}, (r+s)\sgn(\sin2\gamma)\sqrt{|\sin2\gamma|},r-s).
    \end{align*}
    Therefore
    \begin{align*}
        \|\alpha(f)-\alpha(h)\|_2^2 &= (r+s)^2(|\cos2\theta|+|\sin2\theta|) + (r-s)^2\\
        &\geq (r+s)^2 + (r-s)^2 = 2(r^2+s^2).
    \end{align*}
    On the other hand
        $$\min_{|\theta|=1}\|f-\theta h\|=\sqrt{\|f\|^2+\|h\|^2} = \sqrt{r^2+s^2}.$$
    Hence
        $$\min_{|\theta|=1}\|f-\theta h\|\leq \frac{1}{\sqrt{2}}\|\alpha(f)-\alpha(h)\|_2.$$
    When $f=0$, or $h=0$, the stability inequality with constant $\frac{1}{\sqrt{2}}$ follows easily. This constant is optimal, as one can observe by taking $f=(1,0)$ and $h=(0,1)$.
\end{example}
\begin{remark}
    The same example works in the case of $\ell_p$ norm for any $p\in [2,4]$. A similar calculation gives $C_\perp = 2^{-1/p}$. On the other hand, by taking a pair $(f,h)=(u_{\pi/8},u_{\pi/8+t})$ with $t>0$ sufficiently small, we get that $C_{glob}\geq C_{(f,h)}>2^{-1/p}$. 
\end{remark}

\section{Unsigned maps}
\label{sec:unsigned_maps}

We now consider measurements for which the phase of each coordinate is also lost. Thus, instead of $\alpha$ and $\beta$, we observe their coordinatewise moduli $|\alpha|$ and $|\beta|$. We investigate uniqueness and stability for these unsigned maps in both finite and infinite-dimensional Hilbert spaces.

\begin{definition}
    Let $\{A_j\}_{j\in J}\subseteq B(H)$ be $p$-Bessel. We say:
    \begin{enumerate}
        \item $\{A_j\}_{j\in J}$ does $|\alpha|$-PR if for every $f,h\in H$,
    $$|\alpha(f)|=|\alpha(h)| \implies f=\theta h, \text{ for some } |\theta|=1.$$
        \item $\{A_j\}_{j\in J}$ does $|\beta|$-PR if for every $f,h\in H$,
    $$|\beta(f)|=|\beta(h)| \implies f=\theta h, \text{ for some } |\theta|=1.$$
    \end{enumerate}
    Additionally, we say $\{A_j\}_{j\in J}$ does
    \begin{enumerate}
        \item $(p,|\alpha|)$-SPR with constant $C>0$ if
        $$\min_{|\theta|=1}\|f-\theta h\| \leq C\big\||\alpha(f)|-|\alpha(h)|\big\|_p, \qquad f,h\in H,$$
        \item $(p,|\beta|)$-SPR with constant $C>0$ if
        $$\min_{|\theta|=1}\|f-\theta h\|\max_{|\theta|=1}\|f-\theta h\| \leq C\big\||\beta(f)|-|\beta(h)|\big\|_{p/2}, \qquad f,h\in H.$$
    \end{enumerate}
\end{definition}
More explicitly, suppose we are given a system $\{A_j\}_{j\in J}\subseteq B(H)$. Then for any $f\in H$, $|\alpha(f)|$ returns measurements $\sqrt{|\langle A_jf,f\rangle|}$, while $|\beta(f)|$ returns measurements $|\langle A_jf,f\rangle|$. Moreover, since $|\beta_j(f)|=|\alpha_j(f)|^2$, the notions of $|\alpha|$-PR and $|\beta|$-PR are equivalent. 

Next, suppose additionally that each $A_j$ is semidefinite. Then for each $j$, there exists $\sigma_j\in \{-1,1\}$, such that for every $f\in H$ we have 
    $$|\alpha_j(f)|=\sigma_j\alpha_j(f), \qquad |\beta_j(f)|=\sigma_j\beta_j(f).$$
In particular, it is easy to see that this implies that a system of semidefinite operators $\{A_j\}_{j\in J}$ does $\alpha$-PR ($\beta$-PR) if and only if it does $|\alpha|$-PR ($|\beta|$-PR). It is also true that $\{A_j\}_{j\in J}$ does $\alpha$-SPR with constant $C>0$ if and only if it does $|\alpha|$-SPR with $C$, and similarly for the maps $\beta$ and $|\beta|$. 

More generally, given any system $\{A_j\}_{j\in J}$, from the reverse triangle inequality we get that    
    $$\big\||\alpha(f)|-|\alpha(h)|\big\|_p \leq \|\alpha(f)-\alpha(h)\|_p$$
and 
    $$\big\||\beta(f)|-|\beta(h)|\big\|_{p/2} \leq \|\beta(f)-\beta(h)\|_{p/2}.$$
Combining the above observations, we obtain the following.
\begin{lemma}
\label{lem:unsigned_implies_signed}
    Let $\{A_j\}_{j\in J}\subseteq B(H)$ be $p$-Bessel.
    \begin{enumerate}
        \item If $\{A_j\}_{j\in J}$ does $|\alpha|$-PR then it does $\alpha$-PR. The same holds for the map $\beta$ and $|\beta|$.
        \item If $\{A_j\}_{j\in J}$ does $(p,|\alpha|)$-SPR with constant $C>0$, then it does $(p,\alpha)$-SPR with the same constant. The same holds for the map $\beta$ and $|\beta|$.
    \end{enumerate}
    The converse of both statements holds whenever $\{A_j\}_{j\in J}$ consists of semidefinite operators.
\end{lemma}

In Proposition \ref{prop:beta_implies_alpha} we showed that $(p,\beta)$-SPR implies $(p,\alpha)$-SPR. The same also holds for the unsigned maps.

\begin{prop}
\label{prop:unsigned_beta_implies_alpha}
    Let $\{A_j\}_{j\in J}\subseteq B(H)$. Suppose $\{A_j\}_{j\in J}$ satisfies a  $p$-Bessel condition with constant $B$. If it does $(p,|\beta|)$-SPR with constant $C$, then it does $(p,|\alpha|)$-SPR with constant $C\sqrt{2}B^{1/p}$.
\end{prop}
\begin{proof}
    As a preliminary observation, for every $f\in H$ we have
    $$\|\alpha(f)\|^p_p = \sum_{j\in J}|\langle A_jf,f\rangle|^{p/2} \leq B\|f\|^p$$
    thus $\|\alpha(f)\|_p\leq B^{1/p}\|f\|.$ 
    
    Fix $f,h\in H$. Since $|\beta_j(f)|=|\alpha_j(f)|^2$, for every $j$ we have 
    \begin{align*}
    \big||\beta_j(f)|-|\beta_j(h)|\big|&=\big||\alpha_j(f)|-|\alpha_j(h)|\big|\cdot \big||\alpha_j(f)|+|\alpha_j(h)|\big|.
    \end{align*}
    Taking the $\ell_{p/2}$ norm and using Hölder's inequality with $\frac{2}{p}=\frac{1}{p}+\frac{1}{p}$ gives
    \begin{align*}
        \big\||\beta(f)|-|\beta(h)|\big\|_{p/2} & \leq \big\||\alpha(f)|-|\alpha(h)|\big\|_p \cdot (\big\||\alpha(f)| + |\alpha(h)|\big\|_p )\\
        &\leq \big\||\alpha(f)|-|\alpha(h)|\big\|_p \cdot (\big\|\alpha(f)\|_p + \|\alpha(h)\big\|_p )\\
        &\leq \big\||\alpha(f)|-|\alpha(h)|\big\|_p \cdot B^{1/p}(\|f\|+\|h\|)\\
        &\leq \big\||\alpha(f)|-|\alpha(h)|\big\|_p \cdot \sqrt{2}B^{1/p} \max_{|\theta|=1}\|f-\theta h\|,
    \end{align*}
    where in the last inequality we applied \eqref{eq:max_phase_inequality}. Suppose now that $|\beta|$-SPR is satisfied with constant $C$. Then
    \begin{align*}
    \min_{|\theta|=1}\|f-\theta h\|\max_{|\theta|=1}\|f-\theta h\| &\leq C\big\||\beta(f)|-|\beta(h)|\big\|_{p/2} \\
    &\leq C\big\||\alpha(f)|-|\alpha(h)|\big\|_p \cdot \sqrt{2}B^{1/p} \max_{|\theta|=1}\|f-\theta h\|.
    \end{align*}
    If $\max_{|\theta|=1}\|f-\theta h\|=0$, then $f=0=h$, and the $|\alpha|$-stability follows. If $\max_{|\theta|=1}\|f-\theta h\|\neq 0$, reducing both sides gives 
    $$\min_{|\theta|=1}\|f-\theta h\| \leq C\sqrt{2}B^{1/p}\big\||\alpha(f)|-|\alpha(h)|\big\|_p.$$
    \end{proof}

Let $H$ be infinite-dimensional. Using Lemma \ref{lem:unsigned_implies_signed} together with Theorem \ref{thm:alpha_beta_failure_on_strong_Bessel} we can deduce that for $p\neq 4$, no $p$-Bessel system does $(p,|\beta|)$-SPR or $(p,|\alpha|)$-SPR. Moreover, Theorem \ref{thm:alpha_beta_semidefinite_failure} shows that such a failure also holds for all semidefinite $4$-Bessel systems.

In Proposition \ref{prop:beta_alpha_stability_on_rank_1} we obtained a $4$-Bessel system of rank-one operators that does $(4,\beta)$-SPR and $(4,\alpha)$-SPR. Nevertheless, this system fails even $|\beta|$-PR and therefore also $|\alpha|$-PR. The same is true for the systems of self-adjoint operators from Theorem \ref{thm:4_Bessel_stability}.

\begin{example}
    Let $\{E_{mn}\}_{m,n\in \N}$ be defined as in Proposition \ref{prop:beta_alpha_stability_on_rank_1}. We have that  
        $$|\beta_{mn}(f)|=|\langle E_{mn}f,f\rangle|=|\langle f,e_n\rangle|\cdot|\langle f,e_m\rangle|.$$
    Take $f=e_1+e_2$ and $h=e_1-e_2$. Then $f\neq \theta h$ for all $|\theta|=1$. Since 
        $$|\langle e_1+e_2,e_k\rangle|=
            \begin{cases}
                1, & k\in\{1,2\},\\
                0, & k\notin\{1,2\},
            \end{cases}$$
    it follows that
        $$|\beta_{mn}(f)|=
            \begin{cases}
                1, & m,n\in\{1,2\},\\
                0, & \text{otherwise}.
            \end{cases}$$
    Similarly, we get 
    $$|\beta_{mn}(h)|=
            \begin{cases}
                1, & m,n\in\{1,2\},\\
                0, & \text{otherwise}.
            \end{cases}$$
    thus $|\beta_{mn}(f)|=|\beta_{mn}(h)|$ for all $m,n$, despite the fact that $f\neq \theta h$ for all $|\theta|=1$.

    Next, let $\{S_{mn}:m<n\}$ and $\{T_{mn}:m<n\}$ for $m<n$ be defined as in Theorem \ref{thm:4_Bessel_stability}. It is easy to see that 
        $$\langle S_{12}f,f\rangle = \sqrt{2}, \qquad \langle S_{12}h,h\rangle=-\sqrt{2}$$
    and $\langle S_{mn}f,f\rangle=0=\langle S_{mn}h,h\rangle$ otherwise. Since both $f$ and $h$ are real, we also have $\langle T_{mn}f,f\rangle=0=\langle T_{mn}h,h\rangle$ for all $m<n$. Therefore the systems 
        $$\{E_{nn}:n\geq1\}\cup\{S_{mn}:m<n\}$$ 
    and 
        $$\{E_{nn}:n\geq1\}\cup\{S_{mn},T_{mn}:m<n\}$$
    still fail $|\beta|$ and $|\alpha|$-PR. Using the same pair of vectors $f$ and $h$, it is easy to check that the perturbed operators from the second part of Theorem \ref{thm:4_Bessel_stability} also fail $|\beta|$ and $|\alpha|$-PR.
\end{example}

\begin{question}
    Does there exist a $4$-Bessel system of indefinite operators on an infinite-dimensional separable Hilbert space that does $(4,|\alpha|)$-SPR or $(4,|\beta|)$-SPR? 
\end{question}

Next, from Theorems \ref{thm:beta_sufficient_on_ortho_pairs} and \ref{thm:alpha_sufficient_on_ortho_pairs} we know that $(p,\alpha)$-SPR and $(p,\beta)$-SPR can be tested on orthogonal pairs. The same can be therefore deduced about $(p,|\alpha|)$-SPR and $(p,|\beta|)$-SPR provided that the system of operators is semidefinite. As it turns out, even if the system consists of self-adjoint operators, such a reduction may fail.

\begin{example}
\label{ex:beta_SPR_on_ortho_but_not_beta_SPR}
    Let $H=\R^2$ and consider the operators
        $$A_1=\begin{pmatrix}
            1&0\\
            0&0
            \end{pmatrix},
            \qquad
        A_2=\begin{pmatrix}
            0&-1\\
            -1&2
            \end{pmatrix},
            \qquad
        A_3=\begin{pmatrix}
            1&0\\
            0&-2
            \end{pmatrix}.$$
    A direct calculation shows that for $z=(x,y)\in\R^2$, we have
        $$|\beta(z)|=\left(x^2,\,|2y^2-2xy|,\,|x^2-2y^2|\right).$$
    We claim that $\{A_j\}_{j=1}^3$ fails $|\beta|$-PR, and thus $|\alpha|$-PR. Let $f=(1,1)$ and $h=(1,0)$. Then $f\neq \pm h$, while
        $$|\beta(f)|=(1,0,1)=|\beta(h)|.$$
    Since $\{A_j\}_{j=1}^3$ fails both $|\alpha|$-PR and $|\beta|$-PR, it also fails $(p,|\alpha|)$-SPR and $(p,|\beta|)$-SPR.

    On the other hand, we claim that the $(p,|\beta|)$-SPR inequality holds uniformly on all orthogonal pairs. This also implies $(p,|\alpha|)$-SPR, as one can observe by reproducing the proof of Proposition \ref{prop:unsigned_beta_implies_alpha} for orthogonal pairs. We first show that $|\beta|$ separates all non-zero orthogonal pairs. Let $u_\gamma=(\cos\gamma,\sin\gamma)$ for $\gamma\in[0,2\pi]$. Let $f,h\in \R^2$ be orthogonal. Without loss of generality, we may assume that $(f,h)\neq (0,0)$, as the stability inequality is trivially satisfied on $(0,0)$ with any constant $C$. Using the same reduction as in Example \ref{ex:stability_constant_not_achieved_on_ortho}, without loss of generality we write $f=ru_\gamma$ and $h=su_{\gamma+\pi/2}$ for some $\gamma\in[0,2\pi)$ and $r,s\geq 0$. 
    
    Suppose that $|\beta(f)|=|\beta(h)|.$ Note that 
        $$|\beta(f)|=r^2(\cos^2\gamma, 2|\sin^2\gamma-\cos\gamma\sin\gamma|,|\cos^2\gamma-2\sin^2\gamma|),$$
    and
        $$|\beta(h)|=s^2(\sin^2\gamma, 2|\cos^2\gamma+\cos\gamma\sin\gamma|,|\sin^2\gamma-2\cos^2\gamma|),$$
    
    \textbf{Case 1.} $\gamma\in \{0,\pi/2,\pi,3\pi/2\}.$ Then either $|\beta_1(f)|=0$ or $|\beta_1(h)|=0$. This gives $rs=0$. The equality of $|\beta_2(f)|=|\beta_2(h)|$ and $|\beta_3(f)|=|\beta_3(h)|$ implies $r=0$ and $s=0$, and therefore $f=0$ and $h=0$.
    
    \textbf{Case 2.} $\gamma\notin \{0,\pi/2,\pi,3\pi/2\}$. Since
        $$r^2\cos^2\gamma = s^2\sin^2\gamma,$$
    and since $(f,h)\neq (0,0)$, we get that $r,s>0$. Put $R:=\frac{\cos\gamma}{\sin\gamma}.$ From $|\beta_1(f)|=|\beta_1(h)|$ we get $\frac{s^2}{r^2}=R^2.$ Next, from $|\beta_3(f)|=|\beta_3(h)|$ we obtain 
        $$r^2\sin^2\gamma|R^2-2| = s^2\sin^2\gamma|1-2R^2|.$$
    Using $\frac{s^2}{r^2}=R^2$ we obtain
        $$R^2|1-2R^2|=|R^2-2|.$$ 
    Therefore 
        $$(R^2-2)^2=R^4(1-2R^2)^2.$$
    Equivalently,
        $$-4(R-1)(R+1)(R^2+1)(R^4-R^2+1)=0.$$
    Since $R^2+1>0$ and $R^4-R^2+1>0$, this gives $R\in\{-1,1\}.$

    Finally, since $|\beta_2(f)|=|\beta_2(h)|$ we get 
        $$r^2|\sin^2\gamma-\cos\gamma\sin\gamma|=s^2|\cos^2\gamma+\cos\gamma\sin\gamma|.$$
    Dividing by $\sin^2\gamma$ and using the fact that $\frac{s^2}{r^2}=1$ gives
        $$|1-R|=|R^2+R|,$$
    which fails to be true for both $R=1$ and $R=-1$. Thus $|\beta(f)|\neq |\beta(h)|$ for every non-zero orthogonal pair $f,h$.

    Hence, the system does $|\beta|$-PR on orthogonal pairs. Now, consider the set 
        $$K=\{(f,h)\in \R^2\times \R^2: f\perp h, ~\|f\|^2+\|h\|^2=1\}.$$
    Note that $K$ is a compact set in a finite-dimensional space. In particular
        $$c_{p}:=\min_{(f,h)\in K}\||\beta(f)|-|\beta(h)|\|_{p/2}>0.$$
    Moreover, if $f\perp h$, we have 
        $$\min_{|\theta|=1}\|f-\theta h\|\max_{|\theta|=1}\|f-\theta h\| = \|f\|^2+\|h\|^2.$$
    By scaling we get that $\{A_j\}_{j=1}^3$ satisfies $(p,|\beta|)$-SPR on orthogonal pairs with constant $1/c_{p}$. Nevertheless, $\{A_j\}_{j=1}^3$ fails $(p,|\beta|)$-SPR.
\end{example}

We know that every system that does $(p,|\beta|)$-SPR also does $(p,|\alpha|)$-SPR. The proof of the converse for unsigned maps uses orthogonal reduction, see the proof of Corollary \ref{cor:beta_and_alpha_are_equivalent}. Since that reduction can fail for unsigned measurements, the converse implication for those maps remains to be considered.

\begin{question}
    Is it possible to find $p\in [2,\infty)$ and an indefinite $p$-Bessel system $\{A_j\}_{j\in J}\subseteq B(H)$ that does $(p,|\alpha|)$-SPR but not $(p,|\beta|)$-SPR?
\end{question}

Finally, from Theorem \ref{thm:beta_and_alpha_PR_in_finite_dim} we know that in finite-dimensional Hilbert spaces, $\alpha$-PR and $\beta$-PR automatically imply $(p,\alpha)$-SPR and $(p,\beta)$-SPR. Therefore, the same implication holds for the maps $|\alpha|$ and $|\beta|$, when every measurement operator is semidefinite. The following example shows that in general stability need not be achieved, even if given a system of self-adjoint operators.

\begin{example} 
\label{ex:beta_pr_not_spr} 
    Let $H=\R^2$ and consider
        $$ A_1= \begin{pmatrix} 1&0\\ 0&1 \end{pmatrix}, \qquad A_2= \begin{pmatrix} 0&0\\ 0&1 \end{pmatrix}, \qquad A_3= \begin{pmatrix} 0&\frac12\\ \frac12&1 \end{pmatrix}. $$ 
    A direct calculation shows that for every $f=(x,y)\in\R^2$ we have 
        $$ |\beta(x,y)| = \left( x^2+y^2,\, y^2,\, |xy+y^2| \right), $$
    and 
        $$ |\alpha(x,y)| = \left( \sqrt{x^2+y^2},\, |y|,\, \sqrt{|xy+y^2|} \right). $$ 
    We first show that the system does $|\beta|$-PR. Since every $(x,y)$ is uniquely determined up to a sign by $(x^2,y^2,xy)$, it is easy to show that whenever $|\beta(x,y)|=|\beta(u,v)|$, we obtain $(u,v)=\pm(x,y).$ Therefore the system does $|\beta|$-PR. 
    Consequently, the system also does $|\alpha|$-PR. 
    
    We now show that it does not do $|\alpha|$-SPR, and hence also $|\beta|$-SPR by Proposition \ref{prop:unsigned_beta_implies_alpha}. For $0<t<1$, let $f_t=(1,t)$ and $h_t=(1,-t).$ 
    Then 
        $$ \min_{|\theta|=1}\|f_t-\theta h_t\| = 2t.$$ 
    On the other hand, 
        $$ |\alpha(f_t)| = \left( \sqrt{1+t^2},\, t,\, \sqrt{t+t^2} \right) $$ 
    and 
        $$|\alpha(h_t)| = \left( \sqrt{1+t^2},\, t,\, \sqrt{t-t^2} \right). $$ 
    Since $|\alpha(f_t)|-|\alpha(h_t)|$ has only one non-zero coordinate, we get 
        $$ \||\alpha(f_t)|-|\alpha(h_t)|\|_p = \sqrt{t+t^2}-\sqrt{t-t^2} = \frac{2t^2}{ \sqrt{t+t^2}+\sqrt{t-t^2} }.$$ 
    In particular, 
        $$ \||\alpha(f_t)|-|\alpha(h_t)|\|_p \sim t^{3/2}, \qquad t\to 0,$$ 
    whereas 
        $$ \min_{|\theta|=1}\|f_t-\theta h_t\| = 2t. $$ 
    Consequently, 
        $$ \frac{ \min_{|\theta|=1}\|f_t-\theta h_t\| }{ \||\alpha(f_t)|-|\alpha(h_t)|\|_p } \longrightarrow\infty, \qquad t\to 0.$$
    Hence the system does not do $|\alpha|$-SPR. From Proposition \ref{prop:unsigned_beta_implies_alpha}, $\{A_j\}_{j=1}^3$ also fails $|\beta|$-SPR.
\end{example}

\begin{question}
    Which finite families of self-adjoint operators that do $|\beta|$-PR also do $(p,|\beta|)$-SPR? And similarly for the map $|\alpha|$?
\end{question}

\subsection*{Usage of Large Language Models}

During the preparation of this manuscript, the authors used ChatGPT (OpenAI) and Claude (Anthropic) to assist with checking mathematical consistency, generating ideas, finding counterexamples and performing numerical calculations. Those tools were also used in the process of literature searching and improving readability of the paper. The authors reviewed and verified the resulting content and take full responsibility for the manuscript, including the correctness of all mathematical statements and proofs.

\subsection*{Acknowledgments}
 DN was partially supported by NSF DMS 2408912, NSF DMS 2522495, NSERC Discovery Grant (2026-04609), a Tier I Canada Research Chair in Mathematics of AI and Machine Learning (CRC-2025-00410), and the Pacific Institute for the Mathematical Sciences (PIMS). YP was partially supported by NSERC DG (2020-04572) and NSERC DAS (2020-00092). TS was partially supported by NSERC DG (2026-04609). ÖY was supported by an NSERC Discovery Grant (22R82411). ÖY also acknowledges support by the Pacific Institute for the Mathematical Sciences (PIMS) and the CNRS – PIMS International Research Laboratory.

\bibliographystyle{plain}
\bibliography{bibfile-SPRforOperators}

\end{document}